\documentclass[12pt]{article}

\usepackage{amssymb,amsmath,latexsym,graphics, supertabular}
\usepackage[enableskew]{youngtab}
\usepackage{amsbsy}
\usepackage{amscd}
\usepackage{verbatim}
\usepackage{amsthm}
\usepackage{times}
\usepackage{mathrsfs}
\usepackage{verbatim}
\usepackage{graphicx,subfigure}
\usepackage[colorlinks]{hyperref}

\usepackage{txfonts}
\usepackage{url}

\usepackage[english]{babel}
\usepackage[T1]{fontenc}
\usepackage[utf8]{inputenc}
\usepackage{tikz}

\usetikzlibrary{backgrounds,fit,decorations.pathreplacing}
\usetikzlibrary{arrows,shapes,positioning}
\usetikzlibrary{calc,through,backgrounds,chains}

\usetikzlibrary{arrows,shapes,snakes,automata,backgrounds,petri}

\newenvironment{pf*}[1]{\proof[#1]}{\endproof}
\newtheorem{Theorem}[equation]{Theorem}

\newtheorem{Corollary}[equation]{Corollary}
\newtheorem{Lemma}[equation]{Lemma}
\newtheorem{Proposition}[equation]{Proposition}

\theoremstyle{definition}
\newtheorem{Definition}[equation]{Definition}
\newtheorem{Example}[equation]{Example}

\theoremstyle{remark}

\newtheorem{Remark}[equation]{Remark}

\numberwithin{equation}{section}
\numberwithin{figure}{section}

\newcommand{\DD}{\mathbb D}
\newcommand{\C}{{\mathbb C}}
\newcommand{\Z}{{\mathbb Z}}

\newcommand{\N}{{\mathbb N}}

\newcommand{\mc}[1]{\mathcal{#1}} % short for mathcal
\newcommand{\mb}[1]{\mathbb{#1}} % short for mathblackboard
\newcommand{\mt}[1]{\text{#1}}

\newcommand{\gr}{\operatorname{gr}}
\newcommand{\lgth}[1]{| #1 |}
\newcommand{\lgthp}[1]{\lgth{#1}'}

\everymath{\displaystyle}
\def\mN{\widetilde{\N}}
\def\M{{\mathbb M}}

\renewcommand{\P}{\varepsilon_{B_n}}

\title{Generators of the Hecke algebra of $(S_{2n},B_n)$.}
\author{K\"ur\c{s}at Aker, Mahir Bilen Can}

\begin{document}

\maketitle

\begin{abstract}
A set of ring generators for the Hecke algebra of the Gel'fand pair $(S_{2n},B_n)$,
where $B_n$ is the hyperoctahedral subgroup of the symmetric group $S_{2n}$ is presented. 
Various corollaries are given. A conjecture of Sho Matsumoto is proven. 
\end{abstract}

\section{\textbf{Introduction.}}\label{sec:intro}

Recall that the {\em Hecke algebra} of a pair of finite subgroups $K\subseteq G$ is the $\Z$-algebra of functions
$$
H_\Z(G,K) := \{ f:G\rightarrow \Z \ | \ f(kgk')=f(g)\ \text{for all}\ f,f'\in K,\ \text{and}\ g\in G\},
$$
whose multiplicative structure is given by the convolution product.
A case of particular interest, especially for harmonic analysis, is when $K=G^\sigma$ is the subgroup of fixed 
elements of an automorphism $\sigma :G\rightarrow G$ of order 2. 
In this case, $H_\Z(G,K)$ is a finitely generated commutative algebra. See Chapter VII, Example 6 of \cite{Macdonald95}.

Pairs of the form $(G,G^\sigma)$ are called {\em symmetric}, and there are plenty of them. 
For example, for any group $G$, the map $(g,h) \mapsto (h,g)$ is an involution on $G\times G$, 
and furthermore, the fixed subgroup is isomorphic to $G$ in $G\times G$.  Therefore, $(G\times G,G)$ is a symmetric pair.

In this manuscript, we focus on the symmetric pair $(S_{2n},B_n)$, where $G=S_{2n}$, the symmetric group on $[2n]:= \{1,\dots,2n\}$,  
and $K= B_n$, the hyperoctahedral subgroup consisting of permutations $w\in S_{2n}$ that commute with the element 
\begin{align}\label{A:definition of t}
t :=(1\,2)(3\,4)\cdots(2n-1\,2n) \in S_{2n} 	
\end{align}
interchanging $2i-1$ with $2i$ for $i=1,\dots, n$.
We prove that the Hecke ring $H_\Z(S_{2n}, B_n)$ is generated by the elements $H_i$ ($1\leq i \leq n$), sums of $w\in S_{2n}$ 
with exactly $i$ cycles in a graph $\varGamma(w)$ naturally associated with $w$.

Our motivation stems from the classical work of Farahat and Higman \cite{FH59}, and from 
more recent work of Wang \cite{Wang04} on the centers of the integral group rings of wreath products.

According to \cite{FH59}, as a ring, the center of the integral group ring of the symmetric group is generated 
by the elements $Z_i \in \Z[S_n]$, sums of permutations with exactly $i$ cycles (fixed points of permutations 
count as $1$-cycles). 
Indeed, using the embedding of $S_n$ into $S_{n+1}$ as the stabilizer of $n+1$, for $n=1,2,\dots$, 
Farahat and Higman constructs a universal ring $\mc{Z}$ and a filtration on $\mc{Z}$. 
Then, it is shown that the universal ring $\mc{Z}$ is generated by an algebraically independent set of generators 
$T_0, T_1, \ldots$. Existence of surjections from $\mc{Z}$ onto the centers, under which a generator $T_m$ is sent to 
$Z_{n-m}$ for $m=1,2, \ldots$, guarantees that the centers are generated by the $Z_i$'s.

Farahat and Higman's construction of $\mc{Z}$ seems mysterious at first, 
as there are no reasonable homomorphisms between centers of the various group rings of symmetric groups. 
In \cite{MO01}, using Olshanski's centralizer method, Molev and Olshanski provides an alternative, demystifying treatment 
of the construction. In \cite{IK01}, Ivanov and Kerov calculates the structure constants of $\mc{Z}$.

In our work, following the classical approach we construct a universal ring $\mc{H}$ and a filtration on it.
In fact, in Theorem \ref{thm:GrCentHeckeIsom} we observe that the associated graded rings of the universal rings 
$\mc{Z}$ and $\mc{H}$ are isomorphic to each other. 
As a corollary of our result and Theorem 2.5 of \cite{FH59}, we see that 
$\mc{H}$ is free polynomial algebra over countably many indeterminates (Theorem \ref{thm:HeckeIsAFreePolynomialRing}).
Our main result, Theorem \ref{thm:GensForHecke_n} states that the elements $H_i$ generate $H_\Z(S_{2n}, B_n)$ as a ring.

There are some important corollaries of our results. Let us mention a few of them. 
First of all, it is desirable to know more about the generators. 
For $Z_i$'s, a fundamentally important result of Jucys \cite{Jucys} and Murphy \cite{Murphy} states that 
\begin{equation} \label{eqn:JMSn}
Z_i = e_{n-i}(J_1, J_2, \ldots, J_n) ,\ i=1,\dots, n,
\end{equation}
where $e_k$ denotes the $k$-th elementary symmetric function, and $J_1, \ldots, J_n$ are the \emph{Jucys-Murphy elements}, 
defined by $J_1=0$ and $J_k:=(1,k)+(2,k)+\cdots+(k-1,k)$ for $k=2,\ldots,n$.

Let $\varepsilon_{B_n}$ denote the average of the elements of $B_n$. 
From a completely different perspective than ours, the following analog of (\ref{eqn:JMSn}) is observed by Zinn-Justin in \cite{Zinn} 
and Matsumoto in \cite{Matsumoto}:
\begin{equation*}
H_i = e_{n-i}(J_1,J_3,\ldots,J_{2n-1})\cdot \varepsilon_{B_n} = \varepsilon_{B_n} \cdot e_{n-i}(J_1,J_3,\ldots,J_{2n-1}).
\end{equation*}
Therefore, it follows from our Theorem \ref{thm:GensForHecke_n} that Conjecture 9.1 of \cite{Matsumoto} is true, 
namely, the map which sends a symmetric polynomial, $F$, in $n$ variables with integral coefficients to 
$F(J_1,J_3,\ldots, J_{2n-1})\cdot \varepsilon_{B_n}$ is \emph{surjective}. We record this as our Theorem \ref{thm:SymmHecke}. 
In fact, we prove more than just this conjecture; in Theorem \ref{thm:CentHeckeFinalIsom} 
we show that the universal rings $\mc{Z}$ and $\mc{H}$ are isomorphic.

The associated graded ring, $\gr\mc{Z}$ of $\mc{Z}$ is of interest to geometers. 
The ring structure of the total cohomology space $\bigoplus_{n=1}^\infty H^*(Hilb^n(\C^2), \Z)$ 
of the Hilbert scheme of points in the affine place $\C^2$ is studied by Lehn and Sorger in \cite{LS01} and independently, by Vasserot in \cite{Vasserot}.
It is recognized by Wang that this total cohomology ring is isomorphic to $\gr\mc{Z}$.
In \cite{Wang04}, Wang extends the construction of $\mc{Z}$ to the wreath products of symmetric groups with finite groups. 
Moreover, he shows that these universal rings are isomorphic to the cohomology rings of Nakajima quiver varieties, 
as well as to the cohomology rings of certain other Hilbert schemes. 
In \cite{LQW04}, Li, Qin and Wang show that the total cohomology ring for an arbitrary smooth quasi-projective 
complex algebraic surface is a Farahat-Higman type ring, that is to say, isomorphic to the associated graded ring of a universal ring.

As a corollary of our Theorem \ref{thm:GrCentHeckeIsom}, which is mentioned above, we show that 
$\gr\mc{H}$ is isomorphic to the total cohomology ring $\oplus_{n=1}^\infty H^*(Hilb^n(\C^2), \Z)$. 
It would be of interest to pursue further geometric implications of this isomorphism.

We structured our paper as follows. In Section 2 we lay out the notation for symmetric groups and 
introduce concepts in the order which we generalize them to the pair $(S_{2n},B_n)$ in Section 3. 
In Section 4, we use developments of Section 3 to state our results. In Sections 5 and 6 our claims are proved.

\hspace{1cm}

\textbf{Acknowledgements.} We would like to thank Omar Tout for his remarks on an earlier 
version of the manuscript which improved the quality of the paper.

\section{\textbf{Notation and Preliminaries for Symmetric Groups.}}

In this section we set our notation and paraphrase the preliminary steps of \cite{FH59}.

\subsection{Partitions.}

The set of positive integers $\{1,2,3,\dots \}$ is denoted by $\N$, and $[n]$ denotes the finite set $\{1,2,\dots, n\}$.

A {\em partition} is a non-increasing sequence of non-negative integers. We denote by $\mc{P}$ the set of all partitions. 
Given  $\lambda = (\lambda_1,\dots, \lambda_k) \in \mc{P}$, its entries are called {\em parts}, 
and its {\em size} is defined to be the sum of its parts. We denote the size of $\lambda$ by $|\lambda |$. 
The number of parts of $\lambda$ is called its {\em length} and denoted by $\ell(\lambda)$. 
The {\em weight} of $\lambda$ is defined to be the sum $\mt{wt} (\lambda) := | \lambda | + \ell(\lambda)$.

For $\mu, \lambda \in \mc{P}$ we write $\mu \subset \lambda$, if all parts of $\mu$ appear in $\lambda$.
In this case $\lambda - \mu$ denotes the partition obtained from $\lambda$ by removing the parts of $\mu$. 
We use the notation $\mu \cup \lambda$ to denote the partition obtained by first juxtaposing $\lambda$ and $\mu$ together and 
then ordering the result in a non-increasing fashion. 
The partition $\lambda + \mu $ is the partition obtained by vector addition. 
For example, $(3,3,1) \cup (5,2) = (5,3,3,2,1)$ and $(3,3,1) + (5,2) = (8,5,1)$.

For $i \in \N$, let $m_i = m_i (\lambda)$ denote the number of times the integer $i$ appears as a part of $\lambda$. 
Then $\lambda= (1^{m_1} 2^{m_2} 3^{m_3} \cdots)$ is an alternative notation to $\lambda = (\lambda_1,\dots, \lambda_k)$.

For a positive integer $n$ greater than $\mt{wt}(\lambda)$, the {\em $n$-completion} $\lambda(n)$ of $\lambda$ is 
defined by $\lambda(n) := \lambda + (1^{n- |\lambda |})$.
It follows from the definitions that $|\lambda (n) | = n$ and that $\ell(\lambda(n)) = n- |\lambda |$.

\subsection{The infinite symmetric group and stable cycle type.}

We represent permutations of $[n]$ in one-line notation as follows. If $x\in S_n$ maps $i$ to $x(i)$ for $i=1,\dots,n$, then 
we write $x=(x_1, x_2, \dots, x_n)$. We multiply permutations from right to left. For example, if $x = (2,1,3)$ and $y =(3,2,1)$, 
then $x y =(3,1,2)$.

A string of distinct positive integers $(i_1 i_2 \cdots i_r)$ is called a {\em cycle} of length $r$. 
Cycles $C_1$ and $C_2$ are called {\em disjoint}, if they have no common entries.
It is well known that every permutation $x\in S_n$ is a concatenation $C_1 C_2 \cdots C_k$ of disjoint cycles, 
for which every number $i\in [n]$ appears in exactly one of $C_1, \dots, C_k$, and if the number $i$ is followed by $j$ 
in a given cycle $C_s$, then $x (i) = j$.

The {\em cycle type} of an element $x =C_1C_2 \cdots C_k \in S_n$ is defined to be the partition obtained by reordering 
the lengths of $C_i$'s in a non-increasing fashion. 
Two permutations $x,y\in S_n$ are called {\em conjugate}, if there exists $g\in S_n$ such that $y= gxg^{-1}$. 
Conjugation preserves the cycle type.

Let $S_\infty$ denote the {\em infinite symmetric group}, which, by definition, is the group of finite-support automorphisms of $\N$. 
In other words, a bijection $x : \N \rightarrow \N$ lies in $S_\infty$ if and only if $x(i) = i$ for all but finitely many $i\in \N$. 
Therefore, each $x \in S_\infty$ lies in some $S_n$.

The {\em stable cycle type} of a permutation $x$ is defined to be the partition 
$$
\lambda -( 1^{\ell(\lambda)}) = (\lambda_1 -1, \lambda_2 -1,\dots, \lambda_{\ell(\lambda)} -1),
$$
where $\lambda$ is the cycle type of $x$ viewed as an element of $S_n$ for some $n \geq 1$. 
Note that the stable cycle type is independent of the symmetric group that $x$ sits in.

\subsection{The Cayley graph and its degree function.}

Let $G$ denote a group and $T\subseteq G$ a generating set. 
The {\em Cayley graph} of the pair $(G,T)$ is the directed graph $\mc{C}(G,T)$ with the vertex set $G$ and 
edges defined as follows. Let $x$ and $y$ be two elements from $G$. Then there exists a directed edge from $x$
to $y$, if $yx^{-1} \in T$. In this case, we write $x\to y$.

The {\em distance function} $d : G \times G \rightarrow \N \cup \{0\}$ is defined so that, for $x,y\in G$,  
$d(x,y)$ equals the length of a shortest directed path from $x$ to $y$. 
Let $e\in G$ denote the identity element. The {\em degree} of $x\in G$, denoted by $|x|$, is defined to be   
$$
|x|:=d(e,x).
$$	
Notice that $|x|$ is equal to the number of elements that appear in a {\em minimal decomposition} $x=t_1 t_2 \cdots t_m$, 
where $t_i \in T$ for $i=1,\dots, m$. Hence, $m=|x|$.

\begin{Example}
Let $T \subseteq S_n$ denote the set of all transpositions. Clearly $T$ generates $S_n$. 
We have the following observations for the Cayley graph of the pair $(S_n,T)$:
\begin{enumerate}
\item The degree of $x \in S_n$ is the minimal number of transpositions whose product is equal to $x$. 
\item If $\lambda$ is the stable cycle type  of $x$, then $|x|= |\lambda|$.
\item For all $z\in S_n$, $|x| = | zxz^{-1}|$.
\item For any $x,y\in S_n$, $|xy| \leq |x| + |y|$.
\end{enumerate}
\end{Example}

The proof of the next proposition is standard, so we omit it. 
\begin{Proposition} \label{prop:Cayley1} 
Let $T^{-1}$ denote the set of elements $t^{-1} \in G$, where $t\in T$. Then 
\hskip 1pt

\begin{enumerate}

\item The right action of $G$ on itself induces a distance preserving action on $\mc{C}(G,T)$.
	
\item If $T=T^{-1}$, then $x \to y$ if and only if $y \to x$. Furthermore, in this case, 
	\begin{enumerate}
	\item $\mc{C}(G,T)$ is undirected,
	\item $d(x,y)=d(y,x)$ for all $x,y \in G$,
	\item $|x| = |x^{-1}|$ for all $x\in G$,
	\item $|xy|\leq |x| + |y|$ for all $x,y \in G$.
	\end{enumerate}

\item If $T$ is invariant under conjugation, then the left action of $G$ on itself induces a left action on $\mc{C}(G,T)$. 
In this case, 	
	\begin{enumerate}
	\item the action is distance-preserving, 
	\item the degree of two conjugate elements are the same. Equivalently,
	the degree $|\cdot|$ is constant along the conjugacy classes.
	\end{enumerate}

\end{enumerate}
\end{Proposition}

\begin{comment}
\begin{proof}
(1) Given an edge $x \to y$ of $\mc{C}(G,T)$ and an element $g \in G$, 
there exists an edge $xg  \to yg$, because $yg(xg)^{-1} = yx^{-1} \in T$. 
This proves the right action of $G$ induces an action on the set of edges of
$\mc{C}(G,T)$ compatible with the action on the vertices. Any group action on a graph 
is distance preserving.

(2) Recall that there is an edge $x \to y$ in the Cayley graph if and only if $yx^{-1} \in T$.
Since $T=T^{-1}$, the element $xy^{-1} = (yx^{-1})^{-1} \in T$. Hence, there is an
edge $y \to x$ in the Cayley graph.

(2a-b) The graph is undirected and $d(x,y)=d(y,x)$
for all $x,y \in G$.

(2c) For any $x\in G$,
$|x| = d(e,x) = d(x^{-1}, e) = d(e, x^{-1}) = |x^{-1}|.$

(2d) For any $x,y \in G$,
$$
|xy|  = d(e,xy) = d( y^{-1}, x ) \leq d(y^{-1}, e ) + d(e,x) = |x| + |y|.
$$ 

(3) Let $x \to y$ be an edge, hence, $yx^{-1} \in T$. By hypothesis, $T$ is preserved under conjugation. 
Then for any $g \in G$, $(gy)(gx)^{-1} = g(yx^{-1})g^{-1}\in T$. Therefore, there is an edge $gx \to gy$. 
In other words, $G$ has a natural left action on $\mc{C}(G,T)$.

(3a) Distance preservation is automatic. 

(3b)
For any $g, x\in G$, 
$$
|g^{-1} x g| = d(e,  g^{-1} x g ) = d(g, xg) = d( e, x ) = |x|.
$$

\end{proof}
\end{comment}

Let $G$ be the infinite symmetric group $S_\infty$ and let $T$ be the set of all transpositions $\{ (i,j) : i \neq j \}$. 
Then, $T=T^{-1}$. Moreover, $T$ is the conjugacy class of elements whose stable cycle type is $(1)$.
\begin{Lemma}
If $x \in S_\infty$ has the stable cycle type $\lambda$, then $|x|=|\lambda|$.
\end{Lemma}
\begin{proof}
Let $x\in S_\infty$ be an element of stable cycle type $\lambda$.
It suffices to check the equality for cycles only. It is well known that any $(r+1)$-cycle $x=(i_0\, i_1 \ldots i_r)$ 
can be written minimally as the product of $r$-transpositions, $x=(i_0\, i_r) \cdots  (i_0\, i_2) (i_0\, i_1)$. 
Therefore, $|x| = r$. On the other hand, since the stable cycle type  of $x$ is $(r)$, we see that $|x| = r = |(r)|$, as claimed.
\end{proof}

\subsection{The support of a permutation.}

The (ordinary) {\em support} of a permutation $x\in S_n$ is defined by $\N(x) := \{ i\in \N:\ x(i) \neq i \}$.
The proof of the following lemma is omitted.

\begin{Lemma}\label{R: stable cycle type vs support}
Let $\mu$ denote the stable cycle type of $x$. Then 
\begin{itemize}
\item $|\N(x) | = \mt{wt}(\mu)= |\mu| + \ell(\mu)$. 
\item For $x, z \in S_{\infty}$, $\N(zxz^{-1}) = z\N(x)$. 
\item If $\N(x_1) \subseteq \N(x_2)$ for some $x_1, x_2 \in S_{\infty}$, then $\N(zx_1z^{-1}) \subseteq \N(zx_2z^{-1})$.
\end{itemize}
\end{Lemma}

The {\em support of an $r$-tuple} $(x_1, \ldots, x_r)\in S_n \times \cdots \times S_n$ is the union of supports of its entries: 
$\N (x_1, \ldots, x_r ) := \bigcup_{i=1}^r\, \N(x_i)$.
We call $r$-tuples $(x_1, \ldots, x_r)$ and $(y_1,\ldots, y_r)$ {\em simultaneously conjugate}, 
if there exists a single element $z$ such that $x_i = z y_i z^{-1} $ for all $i=1, \ldots, r$. 
It is immediate from Lemma \ref{R: stable cycle type vs support} that the size of $\N(x_1, \ldots, x_r )$ 
is invariant under simultaneous conjugation.

\begin{Lemma} \label{lem:Support} \label{lem:movers}
For $x,y \in S_{\infty}$,
\begin{enumerate}
\item[(1)] $|\N(xy)| \leq |\N(x,y^{-1})|$.
\item[(2)] In (1), the equality holds if and only if $\N(y) \subseteq \N(xy)$.
\end{enumerate}
\end{Lemma}

\begin{proof}
Denote the complement of a subset $X\subseteq Y$ by $X^c$.
(1) Let $y^{-1} \sqcup 1$ denote the map
\begin{align*}
y^{-1} \sqcup 1:\ \N(y^{-1})^c \longrightarrow \N(xy)^c \sqcup ( \N(x) \setminus \N(y^{-1}) )
\end{align*}
defined as follows. 
If $i \in \N(y^{-1})^c$, then either $i \in \N(x)$, or not. If $i \in \N(x)$, then $i \in \N(x) \setminus \N(y^{-1})$ 
and the maps sends to $i$ to itself. If not, then $i \in \N(x)^c$, and $y^{-1} i \in \N(y)^c$, therefore 
$y^{-1} i \in \N(xy)^c$. In this case, the map sends $i$ to $y^{-1} i$.

Note that $\N(x) \setminus \N(y^{-1}) \subseteq \N(xy)$. 
Therefore, the target is a disjoint union of two sets. Note also that, by construction, $\N(x) \setminus \N(y^{-1})$ 
is always in the image. 
Finally, note that $y^{-1} \sqcup 1$ is injective by definition.

In order to obtain the inequality, we assume that $x,y\in S_{n}$ for some $n\in \N$. At this level, 
the injectivity of $y^{-1} \sqcup 1$ implies that 
$$
n - |\N(y^{-1})| \leq (n - |\N(xy)|) + |\N(x) \setminus \N(y^{-1})|.
$$
Therefore, 
$|\N(xy)| \leq |\N(y^{-1})| +  |\N(x) \setminus \N(y^{-1})| = | \N(x) \cup \N(y^{-1})| = |\N(x, y^{-1})|$.

(2) Suppose that $\N(y^{-1}) = \N(y) \subseteq \N(xy)$. Note that $\N(x,y^{-1}) = \N(y^{-1}) \sqcup ( \N(x) \setminus \N(y^{-1}))$. 
Let $i \in \N(x,y^{-1})$. If $i \in \N(y^{-1})$, then there is nothing to prove. 
If $i \in  \N(x) \setminus \N(y^{-1})$, then it follows from the definitions that $i\in \N(xy)$. 
Hence, $\N(x,y^{-1}) \subseteq \N(xy)$. Combined with part (1), we get the equality $|\N(xy)| = |\N(x,y^{-1})|$.

Conversely, suppose that $|\N(xy)| = |\N(x,y^{-1})|$. It follows from the proof of part (1) that $y^{-1} \sqcup 1$ is onto. 
Furthermore, in this case, the set $\N(xy)^c$ is equal to 
$\{ y^{-1} i :\ i\in \N(y^{-1})^c\ \text{and}\ i\in \N(x)^c \}$. It follows that  
$\N(xy)  = \N(y^{-1}) \cup \N(x)$, and hence $\N(y) \subseteq \N(xy)$.

\end{proof}

\begin{Lemma}
Let $x, y\in S_\infty$ be two permutations of stable cycle types $\lambda$ and $\mu$, respectively.
Let $\nu$ denote the stable cycle type of $xy$. Then, 
$|\nu| \leq |\lambda | + |\mu |$. Moreover, the equality holds if and only if $\N(y) \subseteq \N(xy)$.
\end{Lemma}
	
\begin{proof}
By Proposition \ref{prop:Cayley1}, $|xy| \leq |x| + |y|$. Therefore, $|\nu| \leq |\lambda | + |\mu |$. 
The equality follows from Lemma \ref{lem:Support}, part (2).
\end{proof}

\subsection{Center of the group ring.}

In this section we review some basic facts about the center of the integral group ring of the symmetric group. 
We denote by $\mc{Z}_n$ the center of $\Z[S_n]$. 

Let $\mu$ denote a partition, and $n$ denote a positive integer. The {\em class sum of type $\mu$}, denoted by $C_\mu(n)$, 
is the sum of permutations $w\in S_n$ with stable cycle type $\mu$. 
Then $C_\mu(n)$ is non-zero if and only if $\mt{wt}(\mu)=|\mu| +\ell(\mu) \leq n$.
Moreover, if $n$ is fixed, the elements $C_\mu(n)$ with $\mt{wt}(\mu) \leq n$ form a $\Z$-linear basis of the center $\mc{Z}_n$.

Let us denote by $a_{\lambda \mu}^\nu(n)$ the structure constants of the class sum basis, defined by the equations: 
$$
C_\mu(n) C_\lambda(n) = \sum_\nu a_{\lambda \mu}^\nu(n) C_\nu (n).
$$
It is shown in \cite{FH59} that 
\begin{itemize}
\item $a_{\lambda \mu}^\nu(n)=0$ if $|\nu | > |\mu |+|\lambda|$, 
\item $a_{\lambda \mu}^\nu(n)$ is independent of $n$ if $|\nu | = |\mu |+|\lambda|$,
\item $a_{\lambda \mu}^\nu(n)$ is a rational function of $n$ if $|\nu | <  |\mu |+|\lambda|$.
\end{itemize}
Furthermore, the main result of \cite{FH59} is that $\mc{Z}_n$ is generated as a $\Z$-algebra 
by the elements
$$
Z_i = \sum_{\mu:\ \ell(\mu)=i} C_\mu (n),\qquad i=1,\dots, n.
$$
Our main results are the analogs of these observations for the Hecke ring of $(S_{2n}, B_n)$ .

\section{\textbf{Invariants of the symmetric pair $(S_{2\infty}, B_\infty)$.}}

The symmetric groups $S_{2n}$, $n=1,2,\dots$  form an inductive family through compositions 
$$
S_{2n} \hookrightarrow S_{2n+1} \hookrightarrow S_{2n+2}.
$$
The inductive limit of the resulting embeddings $S_{2n} \hookrightarrow S_{2n+2}$ is isomorphic to $S_\infty$. 
However, we write $S_{2\infty}$ for the limit to remind us of the particular inductive structure. 
The hyperoctahedral subgroups $B_n \subset S_{2n}$ form an inductive family compatible with the inclusions 
$S_{2n} \hookrightarrow S_{2n+2}$. We denote by $B_{\infty}$ the resulting limit of $B_n$'s.

\subsection{Stable coset types.}\label{subsection: stable coset types}
	
For each permutation $w \in S_{2n}$, we define a graph $\varGamma(w)$ as follows. 
Vertices of $\varGamma(w)$ are located on a circle and each vertex is labelled by a pair of numbers, 
which we call the {\em exterior} and {\em interior} labels. Exterior labels run through natural numbers from $1$ to $2n$ in a clockwise 
fashion around the circle. The interior label for a vertex with exterior label $i$ is $w(i)$. Vertices with partnering exterior 
(resp. interior) labels are linked by an undirected edge.

\begin{figure}[htp]
\begin{center}
\begin{tikzpicture}[scale=.6]

\begin{scope}
\def \n {12}
\def \radius {3.25cm}
\def \margin {8} 
\foreach \s in {1,...,\n}
{
 \node[blue] at ({360/\n * (\s - 1)}:\radius) {$\bullet$};
 \draw[blue] ({360/\n * (2*\s - 1)}:\radius) arc ({360/\n * (2*\s -
1)}:{360/\n * (2*\s)}:\radius);
}
\end{scope}

\begin{scope}
\draw[thick, blue ,-] (0 : 3.25cm) .. controls +(left: 3cm) and
+(down:.1cm) .. (-60: 3.25cm) ;
\draw[thick, blue ,-] (-210 : 3.25cm) .. controls +(right: 3cm) and
+(down:.1cm) .. (-150: 3.25cm) ;
\draw[thick, blue ,-] (-270 : 3.25cm) to (-120: 3.25cm) ;
\draw[thick, blue ,-] (-30 : 3.25cm) .. controls +(left: 2.3cm) and
+(down:.1cm) .. (-90: 3.25cm) ;
\draw[thick, blue ,-] (-240 : 3.25cm) .. controls +(down:.1cm) and +
(right: 3cm)  .. (-180: 3.25cm) ;
\draw[thick, blue ,-] (30 : 3.25cm) .. controls +(left: 0cm) and +
(down: 3.4cm)  .. (60: 3.25cm) ;
\end{scope}

\begin{scope}
\def \n {12}
\def \radius {3.6cm}
\node[blue] at (0:\radius) {$1$};
\node[blue] at (-360/\n:\radius) {$2$};
\node[blue] at (-360/\n*2:\radius) {$3$};
\node[blue] at (-360/\n*3:\radius) {$4$};
\node[blue] at (-360/\n*4:\radius) {$5$};
\node[blue] at (-360/\n*5:\radius) {$6$};
\node[blue] at (-360/\n*6:\radius) {$7$};
\node[blue] at (-360/\n*7:\radius) {$8$};
\node[blue] at (-360/\n*8:\radius) {$9$};
\node[blue] at (-360/\n*9:\radius) {$10$};
\node[blue] at (-360/\n*10:\radius) {$11$};
\node[blue] at (-360/\n*11:\radius) {$12$};
\end{scope}

\begin{scope}
\def \n {12}
\def \radius {2.9cm}
\node[black] at (0:\radius) {$1$};
\node[black] at (-360/\n:\radius) {$7$};
\node[black] at (-360/\n*2:\radius) {$2$};
\node[black] at (-360/\n*3:\radius) {$8$};
\node[black] at (-360/\n*4:\radius) {$6$};
\node[black] at (-360/\n*5:\radius) {$4$};
\node[black] at (-360/\n*6:\radius) {$10$};
\node[black] at (-360/\n*7:\radius) {$3$};
\node[black] at (-360/\n*8:\radius) {$9$};
\node[black] at (-360/\n*9:\radius) {$5$};
\node[black] at (-360/\n*10:\radius) {$12$};
\node[black] at (-360/\n*11:\radius) {$11$};
\end{scope}

\end{tikzpicture}
\label{G:graphofw}
\caption{Graph of $w$.}
\end{center}
\end{figure}
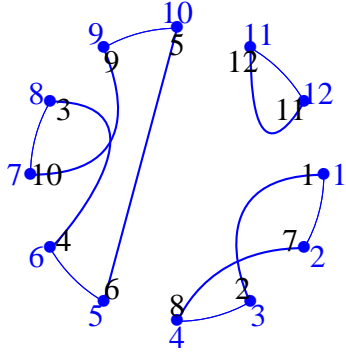

For example, if $w = ( 1, 7, 2, 8, 6, 4, 10, 3, 9, 5, 12, 11)$, then the vertex with exterior label $2$ has 7 as its interior label. 
We depict the graph $\varGamma (w)$ of $w$ in Figure \ref{G:graphofw}.

The graph $\varGamma (w )$ consists of cycles, all of which have even length, say
$2\rho_1 \geq  2\rho_2 \geq \ldots$.
Then the sequence $(2\rho_1, 2\rho_2, 2 \rho_3, \ldots)$ is a partition of $2n$ denoted by $2\rho(w)$, 
where $\rho (w) = (\rho_1,\rho_2,\ldots )$ is a partition of $n$. The partition $\rho(w)$ is called the {\em  coset type of $w$}, because
	
\begin{Lemma} \label{lem:CosetType}
Two elements $x$ and $y$ of $S_{2n}$ are in the same $B_n$-double coset if and only if they have the same coset type.
\end{Lemma}
\begin{proof}
See [VII, \S 2.1, \cite{Macdonald95}].
\end{proof}

Let $w\in S_{2\infty}$ be a permutation. Then $w\in S_{2n}$ for some $n\in \N$. The coset type $\rho(w)$ of $w$ in $S_{2n}$ is well-defined, 
but depends on the integer $n$. To remove this dependency, just as in the case of stable cycle type, we define 
the {\em stable coset type} of $w$ to be 
$$
\overline{\rho} (w):= \rho(w) -( 1^{\ell(\rho (w))}). 
$$
\begin{Example}
Let $z= (2\, 3) (4\,5) (6\, 7) (8\, 9)( 10\,11)(12\,1) \in S_{12}$. Then the stable coset type of $z$ is $(2,2)$.
\end{Example}

There is an obvious generalization of the previous lemma to $S_{2\infty}$.
\begin{Lemma}
Two elements $x$ and $y$ of $S_{2\infty}$ are in the same $B_{\infty}$-double coset if and only if 
they have the same stable coset type .
\end{Lemma}
\begin{proof}
Let $n$ be sufficiently large so that $x,y\in S_{2n}$. Then $x$ and $y$ share the same coset type 
if and only if they are in the same double coset of $B_n$. 
\end{proof}

\subsection{Twisted involution of a symmetric space.}

The involution (order 2 automorphism) of $(S_{2n},B_n)$ is given by the conjugation $\sigma(w) = twt^{-1}$, 
where $t$ is as in (\ref{A:definition of t}).
It is useful to consider a ``twist'' of $\sigma$. More generally, if $(G,K)$ denotes a symmetric pair with $K=G^\sigma$ 
for some involution $\sigma$, then we define its {\em twist} to be the map 
\begin{align}\label{A:varphi}
\varphi (x) = \sigma(x^{-1})x,\ x\in G.
\end{align}

\begin{Example}\label{E:varphi}
Let $\varphi$ denote the obvious twist for $(S_{2\infty},B_\infty)$, and let 
$x$ and $y$ denote the elements $x=(2\, 3) (4\,5) (6\, 7) (8\, 9) \in S_{2\infty}$, 
$y=(2\, 3) (4\,5) (6\, 7) (8\, 9)( 10\,11) \in S_{\infty}$, respectively. Then 
\begin{align*}
\varphi(x) &= \varphi ( (2\, 3) (4\,5) (6\, 7) (8\, 9) (10\,1 )) = (1\, 7\, 3\, 9\, 5 )(2\, 6\, 10\, 4 \,8),\\
\varphi(y) &= \varphi ((2\, 3) (4\,5) (6\, 7) (8\, 9)( 10\,11)(12\,1)) = (1\, 9 \, 5) (11\, 7\, 3)(2\, 6\, 10)(4\, 8\, 12).
\end{align*} 

\end{Example}

Notice that, in general, $\varphi$ is not a group homomorphism. However, the fixed subgroup is a pre-image. 
Indeed, $K= \varphi^{-1}(e)$, where $e$ is the identity element of $G$. 
Given $x,y\in G$, let us write $x\sim y$, if $x$ is conjugate to $y$.

The proof of the following lemma is elementary.

\begin{Lemma} \label{lem:EmbeddingForSymSpace}
For any $x,y\in G$,
\begin{enumerate} 
\item $\varphi(x^{-1})  \sim \varphi(x)^{-1}$,
\item $(\varphi(y^{-1})^{-1}, \varphi(x)\varphi(y^{-1})^{-1}) \sim (\varphi(y), \varphi(xy))$,
\item For all $k\in K$ and $x\in G$, $\varphi(kx)=\varphi(x)$.
\item Right $K$-action on the group $G$ conjugates the image $\varphi(x)$.
\end{enumerate}
\end{Lemma}

\begin{comment}

\begin{proof}
(1) $\varphi(x^{-1}) \sim x^{-1} \varphi(x^{-1}) x = x^{-1} \sigma(x) x^{-1} x = x^{-1} \sigma(x) =\varphi(x)^{-1}$.

(2) Similarly, by conjugating by $\sigma(y^{-1})$ we obtain
$$
\begin{array}{rcl}
(\varphi(y^{-1})^{-1}, \varphi(x)\varphi(y^{-1})^{-1}) & =  & (y\sigma(y^{-1}), \sigma(x^{-1})x y\sigma(y^{-1})) \\
& \sim & (\sigma(y^{-1})y, \sigma(y^{-1}) \sigma(x^{-1})x y) \\
& = & (\varphi(y), \varphi(xy)).
\end{array}
$$

(3) For $k\in K$ and $x\in G$, 
\begin{align*}
\varphi(k x) &=\sigma(x^{-1} k^{-1})(kx) =\sigma(x^{-1})k^{-1} (kx) =\varphi(x).
\end{align*}

(4)
For $h\in K$ and $x\in G$, 
\begin{align*}
\varphi(xh^{-1})=\sigma(hx^{-1})xh^{-1}=h\sigma(x^{-1})  x h^{-1} =h\varphi(x)h^{-1}.
\end{align*}
\end{proof}

\end{comment}

\begin{Remark}
The map $\varphi$  is constant along the left $K$-cosets, hence, it induces a map $K\backslash G \to G$, 
which we denote by $\varphi$, also.
The right $K$-action on the symmetric space $K\backslash G$ alters the image by conjugation. 
Taking quotients, we obtain a sequence of maps,
$$
K\backslash G / K \to G^K \to G^G.
$$
Here, for a subgroup $H$ of $G$, by $G^H$ we denote the set of $H$-conjugacy classes of $G$.

The composition $ K\backslash G / K  \to G^G$ enables us to transfer the invariants of
conjugacy classes of $G$ to the double cosets of $K$. 
When $G=H\times H$ and $K$ is a diagonally embedded copy of $H$ inside $G$, the map 
$K\backslash G / K  \to G^G$ is a bijection.
\end{Remark}

\subsection{Twisted degree.}

Let $(G,K)$ be a symmetric pair and let $T$ denote a conjugation invariant generating subset for $G$. 
We define the {\em twisted degree} of an element $w\in G$ by $|w|' := |\varphi(w) |$, where 
$\varphi$ is the twist defined in (\ref{A:varphi}).

\begin{Proposition} \label{prop:DistanceDegreep}
\hskip 1pt
	
\begin{enumerate}

\item The twisted degree $|\cdot |'$ is constant along the $K$-double cosets.

\item If $x,y\in G$, then $|{x^{-1}}|'= |{x} |'$ and $|{xy}|' \leq |{x}|' + | {y} |'$.
\end{enumerate}

\end{Proposition}

\begin{proof}
(1) Recall that $|\cdot |$ is conjugation invariant. Let $k_1,k_2 $ be two elements from $K$. Then for any $x \in G$ 
$$
| k_1 x k_2 |' = | \varphi(k_1 x k_2) | = |k_2^{-1} \sigma (x^{-1}) x k_2 | = |\varphi(x)| = |x|'. 
$$

(2) For $x\in G$, by Lemma \ref{lem:EmbeddingForSymSpace} (1),
$$ 
|x^{-1} |'= |\varphi(x^{-1})|= |\varphi(x^{-1})^{-1} |= |\varphi(x) | = | x |'.
$$

The inequality $|xy|' \leq |x|' + |y|'$  follows from part (2) of Lemma \ref{lem:EmbeddingForSymSpace} as follows:
$$
|xy|' = |{\varphi(xy)}| = | {\varphi(x) \varphi(y^{-1})^{-1}} | \leq  | {\varphi(x)} | + 
|{\varphi(y^{-1})^{-1}}| = |{x}|' + |{y}|'. 
$$
\end{proof}

We summarize the relationship between the stable coset type and stable cycle type in the following 
\begin{Lemma} \label{lem:CosetTypeVsType}
Let $x\in S_{2\infty}$. Then $x$ is of stable coset type  $\mu$ if and only if the element $\varphi(x)$ is 
of stable cycle type $\mu \cup \mu$. In particular, $|x|'= | {\varphi(x)} | = |{\mu \cup \mu} |= 2 |{\mu} |$.
\end{Lemma}

\begin{proof}
It is enough to prove the claim for single cycles only. To this end, let $y\in S_{2\infty}$ denote an element of stable
coset type $(r)$. For example, let  
$y = (2\,3)\cdots(2r\,2r+1)(2r+2\, 1) \in S_{2r+2}$. 

Divide the set $[2r+2]$ into $4$ subsets $L_i$, $i=0,1,2,3$, with $L_i$ being the set of elements from $[2r+2]$ 
that are equal to $i \pmod{4}$. 

Let $u$ and $v$ denote the two $(r+1)$-cycles defined by:
$$
\begin{array}{rcl}
u & = & (\text{$L_0$ in increasing order}, \text{$L_1$ in decreasing order}) \,\text{and} \\
v & = & (\text{$L_2$ in increasing order}, \text{$L_3$ in decreasing order}).
\end{array}
$$
Then, $\varphi(y)=uv$ and $ |y|' = |{uv} |= |{u} |+ |{v}|=r+r=2r$.
\end{proof}

\begin{Example}
Let $y= (2\, 3) (4\,5) (6\, 7) (8\, 9) (10\,1 ) \in S_{10}$. The stable coset type of $y$ is 
$(4)$ and the stable cycle type of $\varphi(y)= (1\, 7\, 3\, 9\, 5 )(2\, 6\, 10\, 4 \,8)$ is $(4,4)$. 
\end{Example}

\begin{Example}
Let $G= S_{n}\times S_{n}$, $x=(a,b) \in G$ and $\sigma(x)=\sigma(a,b)=(b,a)$. 
Then, $\varphi(x)= \sigma((a,b)^{-1})(a,b)=(b^{-1}, a^{-1})(a,b)=(b^{-1} a, a^{-1} b)$. 
Set $g:=a^{-1} b$. Then, $\varphi(x) = (g^{-1}, g)$. Say $x=( (1,2), e )$, then $|x|'= 2.$
\end{Example}

\subsection{Modified support.}

Consider $X_0$, the set all two-element subsets of $[2n]$. We denote by $X$ the 
set of all $n$ element subsets of $X_0$. There is a permutation action of $S_{2n}$ on $X$.

We call the integers $2j-1, 2j$ {\em partners}, and the pair $\{2j-1, 2j\}$ a {\em couple}. 
The partner of an element $j$ is denoted by $t(j)$.
Let $\mb{D}(n) \in X$ denote the following set of couples:
$$
\mb{D}(n)=\{ \{2j-1, 2j\} \, | \, j \in [n] \}.
$$
Observe that the hyperoctahedral subgroup $B_n$ is the stabilizer of $\mb{D}(n) \in X$ in $S_{2n}$.

Let $\mb{D}$ denote the union of all $\mb{D}(n)$:
$$
\mb{D}:= \{ \{2j-1, 2j\} \, | \, j \in \N \}.
$$
We define the {\em modified support of $x\in S_{2\infty}$} by 
\begin{align}\label{A:support in S_2n,B_n}
{\mN}(x) := \mb{D} \setminus x^{-1}\mb{D} = x^{-1}( x\mb{D} \setminus \mb{D}).
\end{align}
The {\em modified support of an $r$-tuple, $y=(y_1,\ldots, y_r) \in S_{2\infty}^r$} is defined 
to be the union of the modified supports of its entries:
$\mN(Y) := \bigcup_{i=1}^r \mN(y_i)$.

The relation between the modified support (\ref{A:support in S_2n,B_n}) and the ordinary support is as follows. 
\begin{Lemma} \label{L:support vs modified support}
For $x\in S_{2\infty}$, the support $\N(\varphi(x))$ of $\varphi(x)$ is equal to the union of all couples, $C \in \mN(x)$. 
Hence, $|{\N}(\varphi(x))|=2 |{\mN}(x)| $.
\end{Lemma}

\begin{proof}
Let $C = \{ i, t(i)\} $ be a couple. Let us suppose that $xC = \{ x(i), xt(i) \}$ is also a couple.
Then $xt(i)=t(x(i))$. Hence, $(t^{-1} x^{-1} t x)(i)=i$, and vice versa. 
Since $\varphi(x)=t^{-1} x^{-1} t x$, we see that $\varphi(x)$ fixes the point $i$ if and only if 
$x$ sends the couple $\{i, t(i)\}$ onto another couple.
\end{proof}

\subsection{Magnitude of a double coset.}

The group $B_{\infty}\times B_{\infty}$ acts diagonally on $r$-tuples from $S_{2\infty}^r$: 
\begin{align}\label{A:orbit definition}
(a,b) \cdot (z_1, \ldots, z_r) := (az_1b^{-1}, \ldots, az_rb^{-1}).
\end{align}
Obviously, when $r=1$, an orbit of (\ref{A:orbit definition}) is nothing but a $B_\infty \times B_\infty$-orbit.

Since $| {\mN}( z_1, \ldots, z_r  ) |$ is invariant under the action given by (\ref{A:orbit definition}), it makes sense to define 
the {\em magnitude} of an orbit to be the number of elements in the modified support of any of its elements. 
We denote the magnitude of an orbit $L$ by $| \mN (L) |$.

For $n\in \N$, the intersection with  the finite subgroup $S_{2n}^r$ is denoted by $L(n)$: 
$$
L(n) := L \cap (S_{2n} \times \cdots \times S_{2n}).
$$

\begin{Lemma} \label{L:support of a double coset}
Let $x\in S_{2\infty}$ be an element of stable coset type $\mu$ and let $K_\mu$ be its $B_\infty$-double coset. 
Then the magnitude of $K_\mu$ is equal to the weight of $\mu$. In other words, 
$| \mN(x) | = | \mu | + \ell(\mu) = \mt{wt}(\mu)$.
\end{Lemma}
\begin{proof}
By Lemma \ref{L:support vs modified support}, we know that $ | \mN(x) | = \frac{1}{2} | \N(\varphi(x))| $. 
On the one hand,  by Lemma \ref{lem:CosetTypeVsType}, $\varphi(x)$ is of stable coset type $\mu \cup \mu$.
On the other hand, from Lemma \ref{R: stable cycle type vs support} we know that 
the size of the support of an element is equal to the size of its stable coset 
type plus the length of the stable coset type. Therefore, 
$$
| \mN(x) | = \frac{ | \N(\varphi(x))| }{2} = \frac{( |\mu \cup \mu | + \ell(\mu \cup \mu) ) }{2} = |\mu | + \ell(\mu).
$$
\end{proof}

\begin{Lemma} \label{lem: non-empty double cosets}

Let $\mu$ be a stable coset type and let $K_\mu \subset S_{2\infty}$ denote its $B_\infty$-double coset. 
Let $K_\mu(n)$ denote the intersection $K_{\mu}(n) = K_\mu \cap S_{2n}$. 
Then
\begin{enumerate}
\item  $K_{\mu}(n)$ is non-empty if and only if $\mt{wt}(\mu) \leq n$.
\item If non-empty, $K_{\mu}(n)$ is the $B_n$-double coset in $S_{2n}$ corresponding 
to the completion $\mu(n)=\mu + (1^{n-\ell(\mu)})$. 
\end{enumerate}

\end{Lemma}

\begin{proof} 

The second claim follows from definitions and Lemma \ref{lem:CosetType}. We give a detailed 
proof of the first claim. 

($\Rightarrow$) By hypothesis, there exists $w \in K_{\mu} \cap S_{ 2n }$. 
The element $w$ fixes all elements $i>n$ and all couples $\{ 2j-1, 2j \}$ for $j>n$. 
Therefore,  $\N (w) \subset \DD (n)$. Hence, $| \mu | + \ell(\mu)   = \N (w) \leq n$.

($\Leftarrow$) Let $\mu$ be such that $\mt{wt}(\mu)= |\mu | + \ell(\mu)\leq n$ and let $w \in S_{2\infty}$ be an element from $K_\mu$. 
Thus, we know from Lemma \ref{L:support of a double coset} that $|\mN(w)| \leq n$. 
We are going to show that a conjugate of $w$ lies in $K_\mu(n)$. 	
To this end, we define an injection $\delta: S_{\infty} \hookrightarrow B_{\infty}$ as follows.

For $x \in S_{\infty}$, let $x' \in S_{2\infty}$ be the permutation obtained by sending 
$2i-1, 2i$ to $2x(i)-1, 2x(i)$, respectively. Set $\delta(x) = x'$. 
By definition, $\delta(x) \in B_{\infty}$ for all $x\in S_\infty$. Furthermore, by its construction, 
$\delta$ is an injective homomorphism. 

Next, let $X=X_w$ denote the set of $i$ such that $\{ 2i-1, 2i \} \in \mN(w)$. 
By the proof of Lemma \ref{lem:DoubleCosets}, the cardinality of $X$ equals $|\mN (w)|  = |\mu | + \ell(\mu)$. 
Let $u$ denote the permutation such that $ u^{-1}$ maps $X$ bijectively onto $\{1,\ldots, \mN (w) \}$ 
and fixes every other number. Let $v:= \delta (u)$. Then $v^{-1} wv$ fixes all $2j-1, 2j$ point wise, 
except for $j\in \{1,\ldots, |\mN (w)|  \}$, where, by assumption, $| \mN (w)|  \leq n$. 
The element $v^{-1} wv$ fixes all $i\geq 2n$, hence, $v^{-1} wv \in S_{ 2n }$ and the intersection 
$K_{\mu}(n) = K_{\mu} \cap S_{ 2n }$ is non-empty.

\end{proof}

\begin{Lemma} \label{lem:DoubleCosets}

Let $L$ be a $B_\infty \times B_\infty$-orbit as in (\ref{A:orbit definition}) for some $r\geq 1$. 
Suppose that $L(n)$ is a $B_n \times B_n$-orbit in $S_{2n}^{r}$. Then 
\begin{enumerate}
\item[(1)] The number of elements in $L(n)$ is 
$$
\frac{(2^n n!)^2}{k(L) (2^{n-|\mN(L)|} (n- |\mN(L)|)!)}, 
$$
for some nonzero constant $k(L)\in \Z$.

\item[(2)] Suppose that $r=1$ and let $K_\mu $ denote a $B_\infty$-double coset of stable coset type $\mu$. 
Let $K_\mu(n)$ denote the intersection $K_{\mu}(n) = K_\mu \cap S_{2n}$. Suppose also that $K_\mu (n) \neq \emptyset$.
Then the number of elements in $K_{\mu}(n)$ is given by 
\begin{equation} \label{E:numofelts} 
\frac{(2^nn!)^2}{k(K_{\mu}) (2^{n- \mt{wt}(\mu)} (n- \mt{wt}(\mu))!)}, 
\end{equation}	
for some nonzero constant $k(K_{\mu}) \in \Z$.
\end{enumerate}

\end{Lemma}
Note that $L(n)$ might be empty, or it might consist of several $B_n\times B_n$-orbits.

\begin{proof}

(1) By assumption, $L(n)$ is a $B_{n}\times B_{n}$-orbit inside $S_{2n}^{r}$. 
Let $z=(z_1, \ldots, z_r) \in L(n)$ be a point and let $m=|\mN(L)|$
be the magnitude of $L$. Since $z_i$ acts trivially on numbers $j$ when $j > 2m$, the stabilizer 
in $B_n\times B_n$ of $z=(z_1, \ldots, z_r)$ splits:
$$
Stab_{B_{n}\times B_{n}}(z) = Stab_{B_{m} \times B_{m}}(z) \times Stab_{B_{n-m} \times B_{n-m}},
$$
where $B_{n-m}$ stands for the hyperoctahedral group on the symbols $[n] \setminus [m]$.
It follows from definitions that $Stab_{B_{n-m} \times B_{n-m}} \cong B_{n-m}$. Therefore, 
if $k(L)$ denotes the number of elements of the first factor, then the number of elements of the 
orbit $L(n)$ is 
\begin{equation*}
\frac{(2^n n!)^2}{k(L) 2^{n-m} (n- m)!}. 
\end{equation*}
Since $m=|\mN(L)|$, the formula follows.

(2) By Lemma \ref{lem: non-empty double cosets}, when $|\mu | + \ell(\mu)\leq n$, 
$K_{\mu}(n)$ is non-empty and is a single $B_{n}\times B_{n}$-orbit. 
Therefore, the assumptions of part (1) are satisfied. 
The formula follows from part (1) combined with Lemma \ref{L:support of a double coset}.

\end{proof}

Next lemma, which should be compared with [\cite{Wang04}, Lemma 2.8] is about the effect of multiplication on the 
modified support.

\begin{Lemma}\label{L:additivenstat1} 
For any $x,y \in S_{2\infty}$, then $ \lgthp{xy} \leq \lgthp{x} + \lgthp{y}$. 
If the elements $x,y$ and $xy$ are of stable coset types, $\lambda$, $\mu$ and $\nu$, respectively,
then $\lgth{\nu} \leq \lgth{\lambda} + \lgth{\mu}$. Furthermore, 
$\mN (xy) \leq \mN (x,y^{-1} )$ with quality, if $ \lgth{ \nu } = \lgth{ \lambda } + \lgth{ \mu }$.
\label{lem:normIneqB}
\end{Lemma}

\begin{proof} 
First claim follows from Proposition \ref{prop:DistanceDegreep}.
For the second claim, we proceed as follows:
$$
\begin{array}{rcl}
2\mN(xy) & = & \N (\varphi(xy))\\
& = & \N ( \varphi(x)\varphi(y^{-1})^{-1}) \\
& \leq & \N( \varphi(x), \varphi(y^{-1})) \, \text{(by Lemma \ref{lem:movers})} \\
& = & 2 \mN(x, y^{-1}) \, \text{(by definition).}
\end{array}
$$
Therefore, 
$$
\begin{array}{rcl}
\mN(xy) = \mN(x,y^{-1}) & \Leftrightarrow &
\N(\varphi(x)\varphi(y^{-1})^{-1}) = \N( \varphi(x), \varphi(y^{-1})) \\ 
& \Leftrightarrow & \N(\varphi(y^{-1})^{-1}) \subseteq \N(\varphi(x)\varphi(y^{-1})^{-1})\\
& \Leftrightarrow & \N(\varphi(y)) \subseteq \N(\varphi(xy)) \\
& \Leftrightarrow & \N(y) \subseteq \N(xy).
\end{array}
$$ 
\end{proof}

%We record an immediate but intriguing corollary of the previous lemma:
%\begin{Corollary}\label{cor:minimalLength}
%For $x\in S_{2\infty}$, $\lgthp{ x } =   2\min \{ \ell(w):\ w \in B_{\infty} x B_{\infty} \}$.
%\end{Corollary}

\section{\textbf{Main Results.}}\label{S:results}

By abuse of notation, we denote by $K_{\mu}(n)$ the double coset sum corresponding to $K_{\mu}(n)$:
\begin{align*}
K_{\mu} (n) := \sum_{\text{$w$ is of stable coset type $\mu$}} w.
\end{align*}
It follows from Lemmas \ref{lem: non-empty double cosets} and \ref{lem:DoubleCosets} that 
\begin{Corollary}\label{C:main1}
The element $K_{\mu}(n) \neq 0$ if and only if $\mt{wt}(\mu) = |\mu | + \ell(\mu) \leq n$. 
Furthermore, $\{ K_{\mu} (n) : \mt{wt}(\mu) \leq n \}$ is a $\Z$-basis for the Hecke algebra $\mc{H}_n$.
\end{Corollary}

Let $\mc{B}$ be the ring of integer valued polynomials with rational coefficients.
Let $b_{\lambda\mu }^\nu (n)$ denote the structure coefficients defined by the equations 
$K_{\lambda} (n) K_{\mu} (n) = \sum_\nu b_{\lambda\mu}^\nu (n) K_{\nu} (n)$.

\begin{Theorem}\label{T:sabit} 
Structure constants $b_{\lambda\mu}^\nu (n) \in \mc{B}$ of $\mc{H}_n$ obey the tricohotomy:
$$
b_{\lambda\mu}^\nu (n) =
\begin{cases}
0 & ,\text{if}\, \lgth{ \nu } > \lgth{ \lambda } + \lgth{ \mu }, \\
\text{positive integer indep. of $n$} & ,\text{if}\, \lgth{ \nu } = \lgth{ \lambda } + \lgth{ \mu },\\
f^\nu_{\lambda\mu}(n) & ,\text{if}\, \lgth{ \nu } \leq \lgth{ \lambda } + \lgth{ \mu },\\
\end{cases}
$$
where $f^\nu_{\lambda\mu}(t) \in \mc{B}$.
\end{Theorem}

We say that the partitions $\nu \in \mc{P}$ with $ \lgth{ \nu } = \lgth{ \lambda } + \lgth{ \mu }$ form the {\em top part} of the product 
$K_{\lambda}(n) K_{\mu}(n)$, and call the corresponding coefficients $b_{\lambda\mu}^\nu (n)$ as the {\em top coefficients.}

\begin{Definition} 
For partitions $\mu,\nu$ and $\lambda$ from $\mc{P}$ with $\lgth{ \nu } \leq \lgth{ \lambda } + \lgth{ \mu }$,
let $f^\nu_{\lambda\mu}$ denote the constants defined as in Theorem \ref{T:sabit}. 
The {\em universal Hecke ring} $\mc{H}$ is the free $\Z$-module generated by the abstract basis elements $K_\mu$, $\mu \in \mc{P}$, 
satisfying 
$$
K_{\lambda}K_{\mu} := \sum_{ \lgth{ \nu } \leq \lgth{ \lambda } + \lgth{ \mu } } f^\nu_{\lambda\mu} K_{\nu}.
$$
\end{Definition}
By the same reasoning as in [\cite{FH59}, Theorem 2.4], $\mc{H}$ is an associative,  
commutative, unital $\mc{B}$-algebra endowed with natural surjective \emph{ring} homomorphisms 
$\mc{H} \to \mc{H}_n$ for all $n$: 
$$
\sum f_\mu K_\mu \mapsto \sum f_\mu(n) K_{\mu}(n).
$$
The universal Hecke ring is naturally filtered by the above theorem. 
In addition, the associated graded ring $\gr\mc{H}$ has integers for structure constants. 
Following Wang \cite{Wang04}, we call the graded ring $\gr\mc{H}$ the {\em Farahat-Higman ring}.

All of the structure constants for $\gr \mc{H}$ can be calculated once the top coefficients for the product 
$K_{\lambda} (n) K_{\mu} (n)$ is known where $\mu$ has a single part:

\begin{Theorem} \label{thm:ProductExpnForSingleCycle}
Let $r$ denote a nonnegative integer and let $\mu$ and $\nu$ be two partitions satisfying $\lgth{\nu}=\lgth{\mu}+r$. 
Then $K_{\nu}(n)$ appears in the product expansion of $K_{\lambda}(n) K_{(r)}(n)$ 
if and only if there exists a subpartition $\rho$ of $\lambda$ such that $\ell(\rho) \leq r+1$ and 
$ \nu = ( r+\lgth{\rho} )  \cup \lambda - \rho $. 
In this case, the coefficient $b_{\lambda\, (r)}^{\nu}$ is given by
\begin{equation*}
b_{\lambda\, (r)}^{\nu} = \frac{(m_{r+\lgth{ \rho }}(\lambda)+1) (r+\lgth{ \rho }+1 ) r!}{  \prod_{i\geq 0}m_i (\rho)!   },
\end{equation*}
where $m_0(\rho) := r+1 - \ell(\rho)$. In other words, 
\begin{align*}
K_{\lambda}(n) K_{(r)}(n) = \sum_{\rho} b_{\lambda\, (r)}^{( r+\lgth{\rho} )  \cup \lambda - \rho } 
K_{( r+\lgth{\rho} ) \cup \lambda - \rho}(n) ,
\end{align*}
where the sum runs over those partitions $\rho \subseteq \lambda$ so that $\ell(\rho) \leq r+1$.
\end{Theorem}

\begin{Theorem} \label{thm:GrCentHeckeIsom}
The linear map $C_\mu \mapsto K_\mu$ induces an isomorphism of the graded rings $\gr\mc{Z}$ and $\gr\mc{H}$.
\end{Theorem}

\begin{proof}
As in the case of $\gr\mc{H}$, the top coefficients for the product $C_{\lambda} (n) C_{(r)} (n)$ determines the structure constants of 
$\gr\mc{Z}$ in the basis $C_\mu$. By Lemma 3.11, \cite{FH59} or p. 132, \cite{Macdonald95},
these top coefficients are
\begin{align*}
a_{\lambda(r)}^{\nu} &=  \frac{(m_{r+\lgth{ \rho }}(\lambda)+1)  (r+\lgth{ \rho }+1 ) r!    }{  \prod_{i\geq 0} m_i (\rho)! },
\end{align*}
where $\rho$ is a  subpartition of $\lambda$ such that $\ell(\rho) \leq r+1$ and $ \nu = ( r+\lgth{\rho} )  \cup \lambda - \rho $.	
Observe that these coefficients completely agree with $b_{\lambda(r)}^{\nu}$ for all possible $\lambda$, $\nu$ and $r$, 
producing the structure constants for $\gr\mc{Z}$ and $\gr\mc{H}$.
\end{proof}

\begin{Remark}
The linear map $C_\mu \mapsto K_\mu$ does not produce a homomorphism on the level of universal rings $\mc{Z}$ and $\mc{H}$.
To see this, it suffices to compare the products $C_{(1)}C_{(1)}$  and $K_{(1)}K_{(1)}$:
\begin{align*}
C_{(1)}C_{(1)} &=  \frac{t(t-1)}{2} C_{\emptyset} + 3 C_{(2)} + 2 C_{(1,1)},\\
K_{(1)}K_{(1)} &=  t(t-1) K_{\emptyset} + K_{(1)} + 3 K_{(2)} + 2 K_{(1,1)}.
\end{align*}
\end{Remark}

\begin{Remark}
The associated graded ring of $\mc{Z}$ is related to a number of different mathematical objects 
and carries some other natural structures. We name a few. The graded ring $\gr \mc{Z}$ 
\begin{itemize}
\item is isomorphic to $\bigoplus_{n} H^*(Hilb^n(\C^2), \Z)$ \cite{LS01},
\item carries a natural vertex algebra structure (Lascoux-Thibon \cite{LT01}).
\end{itemize}
Thus, by Theorem \ref{thm:GrCentHeckeIsom}, we observe that $\gr \mc{H}$ is isomorphic 
to $\bigoplus_{n} H^*(Hilb^n(\C^2), \Z)$.
\end{Remark}

Recall that the universal ring $\mc{Z}$ is generated by an algebraically independent set of 
generators 
$$
S_i := \sum_{\lgth{\mu}=i } C_{\mu},\  \text{for}\ i \in \N.
$$
See Theorem 2.5 of \cite{FH59}. In our next result, we prove an analog.

\begin{Theorem} \label{thm:HeckeIsAFreePolynomialRing}
The universal Hecke ring $\mc{H}$ is a free polynomial algebra over $\mc{B}$, generated by an algebraically 
independent set of generators $T_i$, $i\in \N$ defined by $T_i := \sum_{\lgth{\mu}=i } K_{\mu}$. 
\end{Theorem}
\begin{proof}
Follows from the proofs of Theorem 2.5 and Lemma 3.15, \cite{FH59}.
\end{proof}

Given a permutation $w\in S_{2n}$, let $c(w)$ denote the number of cycles in the graph $\varGamma(w)$, 
defined in Section \ref{subsection: stable coset types}. We have the following parallel of Theorem 1.1 of \cite{FH59}. 
\begin{Theorem} \label{thm:GensForHecke_n} \label{thm:Main}
The Hecke ring $\mc{H}_n$ is generated by the elements
$$
H_i := \sum_{c(w)=i} w,\ \text{for}\ i=1,\ldots, n.
$$
\end{Theorem}
\begin{proof}
By definition, the element $H_i$ is a linear combination of double coset sums 
$K_{\mu}(n)$ over partitions $\mu$ with $\mt{wt}(\mu) \leq n$:
$$
H_i := \sum_{\lgth{\mu}= n - i} K_{\mu}(n). 
$$

The universal ring $\mc{Z}$ is equipped with a surjective ring homomorphism onto $\mc{H}_n$ which sends
$$
\begin{array}{ccc}
f & \mapsto & f(n), \\
K_{\mu} & \mapsto & K_{\mu}(n), \\
T_i & \mapsto & H_{n-i}.
\end{array}
$$

By the previous theorem, elements $T_i$ generate $\mc{H}$. Therefore, their images $H_{i}$ generate $\mc{H}_n$.
\end{proof}

Let $\Lambda$ denote the ring of symmetric functions with integral coefficients. Our next theorem confirms a 
conjecture of Matsumoto (Conjecture 9.1 in \cite{Matsumoto}).
\begin{Theorem} \label{thm:SymmHecke} 
The ring homomorphism
$$
\begin{array}{ccc}
\Lambda & \longrightarrow & \mc{H}_n \\
F & \mapsto & F(J_1, J_3, \ldots, J_{2n-1}) \cdot \P
\end{array}
$$
is surjective. Under this map, the elementary symmetric polynomial $e_k$ is sent to $H_{n-k}$ for $k\in \{0,\ldots, n-1\}$.

The Hecke ring $\mc{H}_n$ is generated by $e_{k}(J_1, J_3, \ldots, J_{2n-1}) \cdot \P$ for $k=0, \ldots, n-1$,
and for all $n\in \N$.
\end{Theorem}

\begin{Remark}
For any other basis $b_\lambda$ of the ring $\Lambda$, such as monomial symmetric polynomials, 
complete symmetric polynomials, or Schur polynomials, evaluations of elements $b_0, b_1, \ldots, b_{n-1}$ at $(J_1, J_3,\ldots, J_{2n-1})$ 
followed by a product with $\P$ yield a new generating set for $\mc{H}_n$.
\end{Remark}

\begin{Remark}
In \cite{Matsumoto}, effects of the elements of the form $ F(J_1, J_3, \ldots, J_{2n-1}) \cdot \P$ on $\mc{H}_n$ are studied. 
The action of such an element on the spherical function basis is called a spherical expansion. 
Similarly, the action on the double coset basis is called a double coset expansion. 
Our Theorem \ref{thm:SymmHecke} strengthens the results of \cite{Matsumoto} by showing that all elements of 
$\mc{H}_n$ are of the form $ F(J_1, J_3, \ldots, J_{2n-1}) \cdot \P$ for some symmetric function $F \in \Lambda$.
\end{Remark}

Spherical functions of the pairs $(S_n\times S_n, S_n)$ and $(S_{2n},B_n)$ are given by the augmented Schur functions 
$\tilde{s}_\lambda$ and the Zonal polynomials $Z_\lambda$, respectively. These functions 
fit in the family of Jack functions parametrized by a positive real number $\alpha$, for which they correspond to the values 
$\alpha=1$ and $\alpha=2$ respectively. 
It is natural to ask whether the universal rings $\mc{Z}$ and $\mc{H}$ fit in a larger ring with an algebraic 
indeterminate, for which they correspond to the specializations $\alpha=1$ and $\alpha=2$ respectively. 
We address this question as follows. 

\begin{Theorem} \label{thm:CentHeckeFinalIsom}
Let $\Lambda_\mc{B}$ be the ring of symmetric functions over the ring $\mc{B}$. 
Assigning both $S_i \in \mc{Z}$ and $T_i \in \mc{H}$ the degree i, the arrows in the diagram
are isomorphisms of graded rings: 
\begin{figure}[htp]
\begin{center}
\begin{tikzpicture}[scale=.6]
\begin{scope}
\def \n {12}
\def \radius {3.6cm}
\node at (-2,0) {$e_i \in \Lambda_{\mc{B}}$};
\node at (2,2) {$S_i \in \mc{Z}$};
\node at (2,-2) {$T_i \in \mc{H}$};
\draw[thick, ->] (-1.4, -0.4) -- (1.2,-1.6);
\draw[thick, ->] (-1.4, 0.4) -- (1.2,1.6);
\draw[thick, ->] (1.9,1.6) -- (1.9, -1.6);
\end{scope}
\end{tikzpicture}
\end{center}
\end{figure}

The isomorphism $\mc{Z} \to \mc{H}$ induces an isomorphism of the associated graded rings
$\gr \mc{Z}$ and $\gr \mc{H}$ of Theorem \ref{thm:GrCentHeckeIsom}.
\end{Theorem}

\section{\textbf{Proofs.}}\label{S:proofs}

We prepare for the proof of Theorem \ref{T:sabit}. 
To this end, we define the {\em reverted action} of $B_\infty\times B_\infty$ on $S_\infty \times S_\infty$ as follows: 
$$
\begin{array}{lcl}
(a,b) \cdot (x,y) & = & (axb^{-1}, bya^{-1}),
\end{array}
$$
where $x,y \in S_{2\infty}$ and $(a,b) \in B_\infty \times B_\infty$. 

\begin{Definition}
For $(x,y)\in S_\infty \times S_\infty$, define $\M(x,y) := \mN(x,y^{-1})$. Then, $| \M (x,y)| = | \mN (x,y^{-1})|$.
For any $p\in S_{2\infty} \times S_{2\infty}$, the number $| \M(p)| $ is invariant under $B_\infty \times B_\infty$-action.
For an orbit $L$ of the reverted action, as before, we set $L(n):=L\cap (S_{2n}\times S_{2n})$. 
\end{Definition}

\begin{Lemma}\label{L:singleorbit}
Let $L$ be an orbit of the reverted action of $B_\infty \times B_\infty$ on $S_\infty \times S_\infty$.
If  non-empty, $L(n)$ is a single orbit of the reverted action of $B_n\times B_n$ on $S_{2n}\times S_{2n}$. 
\end{Lemma}

\begin{proof}
	
We need to show that if $x,y,z,t\in S_{2n}$ and $x= azb^{-1},\ y=bta^{-1}$ for some $a,b\in B_\infty$, 
then there exist $a_0,b_0\in B_n$ such that $x= a_0 z b_0^{-1}$ and $y=b_0 t a^{-1}$. 
Without loss of generality we assume that $a,b \in B_m$ for some $m > n$. We consider $B_n$ in $B_m$ as the stabilizer of 
$[2m]-[2n]=\{2n+1,2n+2, \ldots, 2m-1,2m\}$.

Let us first look at the cases when the restrictions of $a$ and $b$ to $[2m]-[2n]$ give permutations of the set $[2m]-[2n]$.
In this case, we define $a_0$ and $b_0$ to be the permutations on $[2n]$ by restricting $a$ and $b$, respectively. 
Clearly $x = a_0 z b_0^{-1}$ and $y = b_0 t a_0^{-1}$. 

Next we assume that there exists $2j \in [2m]-[2n]$ such that $b^{-1} (2j) \in [2n]$. 
Before we continue, let us look at an example in detail.
	
\begin{Example}\label{E:detailed}
	
Let 
\begin{align*}
a &= (14, 13, 7, 8, 9, 10, 1, 2, 12, 11, 16, 15, 3, 4, 6, 5) \\
b &= (15, 16, 2, 1, 3, 4, 14, 13, 6, 5, 10, 9, 8, 7, 12, 11)\\
z &= (5, 3, 8, 7, 4, 6, 1, 2, 9, 10, 11, 12, 13, 14, 15, 16),
\end{align*}
hence $a,b\in B_8$ and $z$ can be viewed as an element of $S_{10}$. Define
$$
x =a z b = ( 6, 5, 7, 9, 2, 1, 4, 3, 10, 8, 11, 12, 13, 14, 15, 16).
$$
We would like to show that $a$ and $b$ can be replaced by $a_0 \in B_5$ and $b_0\in B_5$ such that $x= a_0 z b_0$. 
It is enough to choose $a_0, b_0 \in B_8$ such that 
$$
\{a_0 ( i):\ i=11,\ldots, 16 \} = \{b_0(i):\ i=11,\ldots, 16 \}= \{ 11,\ldots, 16\}.
$$ 
	
We start with finding the pairs in $b$ that larger than 10 within $( b_1,\ldots, b_{10} )$, and smaller than or equal to 10 within 
$(b_{11},\ldots, b_{16})$: 
\begin{align*}
\{b_1=15,\ b_2=16 \}, &\ \{b_7= 14,b_8=13\},\\
\{b_{11}=10,b_{12}=9\}, &\ \{b_{13}=8,\ b_{14}=7 \}.
\end{align*}
	
We start modifying $b$ by interchanging  $\{b_1=15,\ b_2=16 \}$ with $ \{b_{13}=8,\ b_{14}=7 \}$:
$\tilde{b} = (8,7, 2, 1, 3, 4, 14, 13, 6, 5, 10, 9, 15,16, 12, 11) \in B_8.$
We modify $a$, as well, so that $x = \tilde{a} z \tilde{b}$. Note that, there is no particular reason to pick the above pairs of $b$ 
to interchange. For every choice of such pairs, $\tilde{a}$  can be appropriately chosen. 
	
To this end we let
$$
\tilde{a}= (5, 6, 7, 8, 9, 10, 1, 2, 12, 11, 16, 15, 3, 4, 13, 14) \in B_8.
$$
In other words, we interchange $\{a_{15}=6,\ a_{16} = 5 \}$ and $\{ a_1=14,\ a_2 = 13 \}$. 
Note the indices: $15= z(b_1),\ 16 = z(b_2),\ 1= z(b_{13}),\ 2=z(b_{14})$. Then, 
$$
\tilde{x} = \tilde{a} z \tilde{b} =  (5, 6, 7, 9, 2, 1, 4, 3, 10, 8, 11, 12, 13, 14, 15, 16).
$$
	
Except the first two entries, $\tilde{x}$ is identical with $x$. However, the transposition $(5,6)$ is in $B_8$, and 
$(5,6) \tilde{x} = x$, and furthermore, $a_0= (5,6)\tilde{a},\ b_0= \tilde{b} \in B_8$ are as we would like them to be. 

\end{Example}

We continue with the proof. The idea is already depicted in the above example. 
Suppose that $x= a z b$ where $x,z \in S_{2n} \hookrightarrow S_{2m}$, $a,b \in B_m$.  
Suppose also that $y = b^{-1} t a^{-1}$ where $y,t \in S_{2n}$. 

Let $\{b_{2i-1}, b_{2i} \}$ be a pair of consecutive entries of $b$ such that $2i \leq 2n$ and $b_{2i} \geq 2n$ 
(since $b\in B_m$, this implies that $b_{2i-1} \geq 2n$). Similarly, let $\{b_{2j},b_{2j-1}\}$ be another pair such that 
$2j > 2n$, and $b_{2j}, b_{2j-1} \leq 2n$. We call such pairs of elements \emph{ misplaced pairs}. 

Observe that $a$ being an element of $B_m$ and the fact that $azb(2j) =2j$ force $\{z (b_{2j}),\ z(b_{2j-1}) \} \in \DD(n)$. 
Let $\{2k-1,2k\} =  \{z (b_{2j}),\ z(b_{2j-1}) \}$. Then, $\{a(2k-1), a(2k) \} = \{ 2j,2j-1\}$. Similarly, let $l> n$ be such that 
$\{b_{2i-1}, b_{2i} \} =  \{2l-1,2l \}  \in  [2m]-[2n]$. Then, $\{a(2l-1), a(2l-1) \} = \{x(2i-1),x(2i)\} \in \DD(n)$.	
	
Define $\tilde{a}$ by interchanging the pairs $\{a(2l-1), a(2l-1) \} $ and $\{a(2k-1), a(2k) \}$ in $a$, 
and define $\tilde{b}$ by interchanging the pairs $\{b_{2j},b_{2j-1}\}$ and $\{b_{2i-1}, b_{2i} \}$. Then, 
	
\begin{itemize}
\item If $s \notin \{ 2i-1,2i,2j-1,2j \}$, then $\tilde{a}z \tilde{b} (s) = azb (s )  = x(s)$.
\item If $s\in \{2i-1,2i\}$, then $\tilde{a}z \tilde{b} (s) \in \tilde{a} z\{ b_{2j-1}, b_{2j} \} = \tilde{a} \{2k-1,2k\}  = \{ x(2i-1),x(2i) \}$.
\item If $s\in \{2j-1,2j\}$, then $\tilde{a}z \tilde{b} (s) \in \tilde{a} z\{ b_{2i-1}, b_{2i} \} = \tilde{a} \{2l-1,2l\}  = \{a(2k-1), a(2k) \} = \{ 2j,2j-1\}$.
\end{itemize}

Therefore, possibly up to a transposition, $\tilde{a} z \tilde{b}$ agrees with $x=azb$. 
We need to show that $\tilde{b}^{-1} t \tilde{a}^{-1}$ agrees with $y=b^{-1} t a^{-1}$ up to the same transposition, also.
By the same token, if $h \notin \{ a(2k-1),a(2k),a(2l-1),a(2l) \}$, then $\tilde{b}^{-1}t \tilde{a}^{-1} (h) = \tilde{b}^{-1}t a^{-1} (h) $. 
If, further, $ta^{-1}(h) \notin \{b_{2j-1} b_{2j}, b_{2i}, b_{2i-1}\}$, then $ \tilde{b}^{-1}t a^{-1} (h) = b^{-1} t a^{-1} (h) = y(h)$.
But $ta^{-1}(h) \in \{b_{2j-1} b_{2j}, b_{2i}, b_{2i-1}\}$ if and only if 
\begin{align*}
h &\in \{a t^{-1} (b_{2j-1}), a t^{-1}( b_{2j}), a t^{-1}(b_{2i}), a t^{-1}(b_{2i-1})\}\\
&= \{a t^{-1} (b_{2j-1}), a t^{-1}( b_{2j}), a (b_{2i}), a (b_{2i-1})\}\\
&= \{a t^{-1} b (2j-1), a t^{-1} b(2j), a (b_{2i}), a (b_{2i-1})\}\\
&= \{2j-1, 2j, a (b_{2i}), a (b_{2i-1})\}.
\end{align*}

The last equality is true because $\{2j-1,2j\} = \{y^{-1}(2j-1),y^{-1}(2j) \} = \{ a t^{-1} b (2j-1), at^{-1} b (2j) \}$. 
In other words, $ta^{-1}(h) \notin \{b_{2j-1} b_{2j}, b_{2i}, b_{2i-1}\}$ if and only if 
$$
a^{-1}(h) \notin \{a^{-1}(2j-1), a^{-1}(2j), b_{2i}, b_{2i-1}\} =\{2k-1,2k,2l-1,2l \} .
$$
Therefore, to show that $\tilde{b}^{-1} t \tilde{a}^{-1} = y$, it is enough to compute $\tilde{b}^{-1}t \tilde{a}^{-1} (h)$ for 
$h\in \{ a(2k-1),a(2k),a(2l-1),a(2l) \}$. To this end, let $h = a(2k)$. 
The other cases $h \in \{ a(2k-1), a(2l),a(2l-1) \}$ are similar. 
\begin{align*}
\tilde{b}^{-1}t \tilde{a}^{-1} (h) = \tilde{b}^{-1}t \tilde{a}^{-1} (a(2k)) \in \{ \tilde{b}^{-1}(2l), \tilde{b}^{-1}(2l-1)\} = \{2j-1,2j \}.
\end{align*}
Since $b^{-1}t a^{-1} (h) = b^{-1}t a^{-1} (a(2k)) \in \{ b^{-1}(2k), b^{-1}(2k-1)\} = \{2j-1,2j \}$, we are done.

This computation shows that, possibly up to a transposition, the number of misplaced pairs in both of 
$\tilde{a}$ and $\tilde{b}$ are one less than that in $a$ and $b$.  
Clearly, we can repeat this argument until we reach at $a_0\in B_n$ and $b_0 \in B_n$ 
with no misplaced pairs, and $a_0 z b_0 = x$, $b_0^{-1} t a_0^{-1} = y$ .

\begin{Example}
We continue with the Example (\ref{E:detailed}). Let 
\begin{align*}
y &= (5, 6, 8, 7, 10, 3, 4, 1, 9, 2, 11, 12, 13, 14, 15, 16),\\
t &= (5, 2, 1, 8, 6, 7, 3, 4, 9, 10, 11, 12, 13, 14, 15, 16),\\
a &= (5, 6, 7, 8, 9, 10, 1, 2, 12, 11, 16, 15, 3, 4, 13, 14),\\
b &= (8, 7, 2, 1, 3, 4, 14, 13, 6, 5, 10, 9, 15, 16, 12, 11),\\
x &= (6, 5, 7, 9, 2, 1, 4, 3, 10, 8, 11, 12, 13, 14, 15, 16),\\
z &= (5, 3, 8, 7, 4, 6, 1, 2, 9, 10, 11, 12, 13, 14, 15, 16).
\end{align*}
Then, $y= b^{-1} t a^{-1}$, and $x = a z b$. 
The misplaced pairs in $b$ are $\{b_7= 14, b_8= 13 \}$ and $\{b_{11}=10,b_{12}=9 \}$. Therefore, 
$\tilde{b} = (8, 7, 2, 1, 3, 4, 10,9, 6, 5, 14, 13, 15, 16, 12, 11)$.
To find $\tilde{a}$ we interchange $\{ a (b_{7})=4, a (b_{8})=3 \}$ and $\{a (b_{11})=11, a(b_{12})=12 \}$:
$$
\tilde{a} = (5, 6, 7, 8, 9, 10, 1, 2, 3,4, 16, 15, 11,12, 13, 14).
$$
Then, 
\begin{align*}
\tilde{a} z \tilde{b} &= (6, 5, 7, 9, 2, 1, 4, 3, 10, 8, 12, 11, 13, 14, 15, 16),\\
\tilde{b}^{-1} t \tilde{a}^{-1} &= (5, 6, 8, 7, 10, 3, 4, 1, 9, 2, 12, 11, 13, 14, 15, 16).
\end{align*}	
Finally, note that $x = (11,12) \tilde{a} z \tilde{b}$ and $y = (11, 12) \tilde{b}^{-1} t \tilde{a}^{-1}$.

\end{Example}

\end{proof}

\begin{Proposition}\label{L:num_of_elts} 
Let $L$ be an $B_\infty \times B_\infty$-orbit for the reverted action. Assume that $L(n)$ is nonempty. 
Then, the number of elements in the intersection $L(n)$ is 
\begin{equation*} \label{E:numofeltsRevOrbit}
\frac{(2^n n!)^2}{k(L) (2^{n-|\M(L)|} (n-|\M(L)|)!)} .
\end{equation*}
for some constant $k(L)$. 
\label{lem:numEltsRevtAction}
\end{Proposition}

\begin{proof}
By the previous lemma, $L(n)$ is a $B_n\times B_n$-orbit for the reverted action.
Pick $(x,y)\in L(n)$. Let $L'$ be the orbit of $(x,y^{-1})$ in the straightforward action. 
The map $(u,v) \to (u,v^{-1})$ sends the orbit $L$ (resp. $L(n)$) of the reverted action to the orbit $L'$ 
(resp. $L'(n)$) of the straightforward action in an equivariant manner.

Plugging $|\M(L)|=|\mN(L')|$ and setting $k(L):=k(L')$ in Lemma \ref{lem:DoubleCosets}, we see that 
the number of elements in $L(n)$ is 
\begin{equation*} 
\frac{(2^nn!)^2}{k(L) (2^{n-|\M(L)|} (n- | \M(L)| )!)},
\end{equation*}
as required.

\end{proof}

Recall Theorem \ref{T:sabit}: 
If $b_{\lambda\mu }^\nu (n)$ are defined by the equations 
$K_{\lambda} (n) K_{\mu} (n) = \sum_\nu b_{\lambda\mu}^\nu (n) K_{\nu} (n)$, then 
$$
b_{\lambda\mu}^\nu (n) =
\begin{cases}
0 & ,\text{if}\, \lgth{ \nu } > \lgth{ \lambda } + \lgth{ \mu }, \\
\text{positive integer indep. of $n$} & ,\text{if}\, \lgth{ \nu } = \lgth{ \lambda } + \lgth{ \mu },\\
f^\nu_{\lambda\mu}(n) & ,\text{if}\, \lgth{ \nu } \leq \lgth{ \lambda } + \lgth{ \mu },\\
\end{cases}
$$
where $f^\nu_{\lambda\mu}(t) \in \mc{B}$.

\begin{proof}[Proof of Theorem \ref{T:sabit}] 

Let $\lambda$, $\mu$ and $\nu$ be the stable coset types as given in the hypothesis. 
By Lemma \ref{L:additivenstat1}, we already know that $b_{\lambda\, \mu}^\nu (n)=0$ if 
$|\nu| > |\lambda| + |\mu|$. To prove the other statements, let $\mc{A}$ denote the set of pairs 
$(x,y) \in S_\infty \times S_\infty$ satisfying $x\in K_\lambda,\ y\in K_\mu,\  xy\in K_\nu$. 
Then $\mc{A}$ is stable under the reverted action of $B_\infty \times B_\infty$. 	
Let  $\mc{A}(n)$ denote the intersection $\mathcal{A} \cap (S_{2n} \times S_{2n})$. 
Hence, $b_{\lambda\, \mu}^\nu (n) = | \mathcal{A}(n) | / |K_\nu (n) |$.

Let $\{A_1,\ldots, A_r\}$ denote the set of orbits of $B_\infty \times B_\infty$ in $\mathcal{A}(n)$. 
Then
\begin{align*}
b_{\lambda\, \mu}^\nu (n) = \frac{| \mathcal{A}(n) | }{ |K_\nu (n) |} = 
\sum_{i=1}^r \frac{ |A_i| }{ |K_\nu (n) | }.
\end{align*}
By (\ref{E:numofelts}) and Proposition \ref{L:num_of_elts}, we see that  
\begin{align}\label{A:complicated_quotient}	
\frac{ | A_i |}{ |K_\nu (n) | } &= \frac{2^{2|\M(A_i)|} (n!)^2}{k(A_i) (n-| \M(A_i)|)!}
\left( \frac{2^{2\mt{wt}(\nu)} (n!)^2} {k(K_{\nu}) (n-\mt{wt}(\nu))!}  \right)^{-1} \\ \notag
&=
\frac{ 2^{2\M(A_i)-\mt{wt}(\nu)} k(K_\nu) (n-\mt{wt}(\nu))!}
{ k(A_i) (n-|\M(A_i)|)!},
\end{align}
for some constants $k(K_{\nu})$ and $k(A_i)$ both of which are independent of $n$. 
In fact, we know from the proof of Lemma \ref{lem:DoubleCosets} that 
$k(K_{\nu})$ and $k(A_i)$ are cardinalities of certain stabilizer subgroups $B'$ and $B''$ respectively.

Let $(x_i,y_i)\in A_i$ be an orbit representative for the reverted action of $B_n \times B_n$. 
Note that $x_i y_i \in K_\nu (n)$.
Then $B'$ is the stabilizer of $(x_i,y_i)$ under reverted action, $B''$ is the stabilizer of $x_i y_i$ 
under the ordinary action of $B_n\times B_n$ on $S_{2n}$. We claim that there is a canonical injection of $B'$ into $B''$.
To prove the claim, let $(a,b) \in Stab_{B_n \times B_n} ((x_i,y_i)) $. Then $a x_i b^{-1} = x_i$ and $b y_i a^{-1} = y_i$. 
Notice that the map $(a,b) \mapsto (a,a) \in Stab_{B_n\times B_n}(x_iy_i)$ is an injective group homomorphism. 
It follows that $k(A_i)$ is a divisor of $k(K_\nu (n))$.
	
By Lemma \ref{lem:normIneqB}, $|\M(A_i)| \geq |\mN(K_{\nu} (n) )| = \mt{wt}(\nu)$ with equality if $\lgth{ \nu } = \lgth{ \lambda } + \lgth{ \mu }$.
Therefore, in the case of equality, (\ref{A:complicated_quotient}) is the integer 
$\frac{k(K_\nu)}{ k(A_i) }$, and therefore so is $b_{\lambda\, \mu}^\nu (n)$. This finishes the proof.

\end{proof}

\section{\textbf{Some explicit computations.}}\label{S:explicit}

In this final section of our paper, we prove Theorem \ref{thm:ProductExpnForSingleCycle},  
which is about the coefficients in the expansion of $K_{\mu}(n) K_{(r)}(n)$. We start by analyzing an example. Let 
\begin{align}\label{A:cycles of x and y}
x &= (7,8,2,9,6,5,12,13,4,1,15,16,14,11,10,3)\\
y &=(3,4,12,1,7,8,10,9,15,11,16,6,13,14,5,2),
\end{align}
then 
$xy = ( 2,9,16,7,12,13,4,1,10,15,3,5,14,11,6,8)$.
It is easy to check that the coset types of $x, y$ and $xy$ are, respectively, $(3,2,1,1,1)$, $(4,1,1,1)$ and $(6,2)$.
Furthermore, if $\lambda$, $\mu$ and $\nu$ denote the modified coset types of $x$, $y$ and $xy$, respectively, 
then $\lgth{ \lambda } + \lgth{ \mu } = \lgth{ \nu }$.

\begin{figure}[htp]
\begin{center}
\begin{tikzpicture}[scale=.47]

\begin{scope}
\node at (-5,0) {$\cdot$};
\node at (5,0) {$=$};
\end{scope}

\begin{scope}[xshift=-10cm]

\begin{scope}
\def \n {16}
\def \radius {3.3cm}
\def \margin {8} 
\foreach \s in {1,...,\n}
{
 \node[blue] at ({360/\n * (\s - 1)}:\radius) {$\bullet$};
 \draw[blue] ({360/\n * (2*\s - 1)}:\radius) arc ({360/\n * (2*\s -
1)}:{360/\n * (2*\s)}:\radius);
}
\end{scope}

\begin{scope}
\draw[thick, blue ,-] (0 : 3.3cm) .. controls +(left: 3cm) and
+(down:.1cm) .. (-22.5: 3.3cm) ;
\draw[thick, blue ,-] (22.5: 3.3cm) .. controls +(left: 3cm) and
+(down: 0cm) .. (-180 : 3.3cm) ;
\draw[thick, blue ,-] (45: 3.3cm) .. controls +(left: 2cm) and +(down:
0cm) .. (-67.5 : 3.3cm) ;
\draw[thick, blue ,-] (-45 : 3.3cm) .. controls +(left: 0cm) and +(up:
0cm) .. (-202.5: 3.3cm) ;
\draw[thick, blue ,-] (-90 : 3.3cm) .. controls +(up: 3cm) and +
(down: 0cm)  .. (-112.5: 3.3cm) ;
\draw[thick, blue ,-] (67.5 : 3.3cm) .. controls +(left:0cm) and +
(down: 0cm)  .. (-135: 3.3cm) ;
\draw[thick, blue ,-] (90 : 3.3cm) .. controls +(down:2cm) and +
(down: 0cm)  .. (-157.5: 3.3cm) ;
\draw[thick, blue ,-] (112.5 : 3.3cm) .. controls +(down:3cm) and +
(down: 0cm)  .. (135.5: 3.3cm) ;
\end{scope}

\begin{scope}
\def \n {16}
\def \radius {4cm}
\node[blue] at (0:\radius) {$1$};
\node[blue] at (-360/\n:\radius) {$2$};
\node[blue] at (-360/\n*2:\radius) {$3$};
\node[blue] at (-360/\n*3:\radius) {$4$};
\node[blue] at (-360/\n*4:\radius) {$5$};
\node[blue] at (-360/\n*5:\radius) {$6$};
\node[blue] at (-360/\n*6:\radius) {$7$};
\node[blue] at (-360/\n*7:\radius) {$8$};
\node[blue] at (-360/\n*8:\radius) {$9$};
\node[blue] at (-360/\n*9:\radius) {$10$};
\node[blue] at (-360/\n*10:\radius) {$11$};
\node[blue] at (-360/\n*11:\radius) {$12$};
\node[blue] at (-360/\n*12:\radius) {$13$};
\node[blue] at (-360/\n*13:\radius) {$14$};
\node[blue] at (-360/\n*14:\radius) {$15$};
\node[blue] at (-360/\n*15:\radius) {$16$};
\end{scope}

\begin{scope}
\def \n {16}
\def \radius {2.9cm}
\node[black] at (0:\radius) {$7$};
\node[black] at (-360/\n:\radius) {$8$};
\node[black] at (-360/\n*2:\radius) {$2$};
\node[black] at (-360/\n*3:\radius) {$9$};
\node[black] at (-360/\n*4:\radius) {$6$};
\node[black] at (-360/\n*5:\radius) {$5$};
\node[black] at (-360/\n*6:\radius) {$12$};
\node[black] at (-360/\n*7:\radius) {$13$};
\node[black] at (-360/\n*8:\radius) {$4$};
\node[black] at (-360/\n*9:\radius) {$1$};
\node[black] at (-360/\n*10:\radius) {$15$};
\node[black] at (-360/\n*11:\radius) {$16$};
\node[black] at (-360/\n*12:\radius) {$14$};
\node[black] at (-360/\n*13:\radius) {$11$};
\node[black] at (-360/\n*14:\radius) {$10$};
\node[black] at (-360/\n*15:\radius) {$3$};
\end{scope}

\end{scope}

\begin{scope}
\begin{scope}
\def \n {16}
\def \radius {3.3cm}
\def \margin {8} 
\foreach \s in {1,...,\n}
{
 \node[blue] at ({360/\n * (\s - 1)}:\radius) {$\bullet$};
 \draw[blue] ({360/\n * (2*\s - 1)}:\radius) arc ({360/\n * (2*\s -
1)}:{360/\n * (2*\s)}:\radius);
}
\end{scope}

\begin{scope}
\draw[thick, blue ,-] (22.5 : 3.3cm) .. controls +(left: 3.5cm) and
+(down: .3cm) .. (292.5: 3.3cm) ;
\draw[thick, blue ,-] (0 : 3.3cm) .. controls +(left: 3cm) and
+(down:.1cm) .. (-22.5: 3.3cm) ;
\draw[thick, blue ,-] (45 : 3.3cm) .. controls + (down: 0.5cm) and +
(down: 4cm) .. (112.5: 3.3cm) ;
\draw[thick, blue ,-] (67.5 : 3.3cm) .. controls + (down: 0cm) and
+(down: 3cm) .. (90: 3.3cm) ;
\draw[thick, blue ,-] (135 : 3.3cm) .. controls +(right: 2cm) and
+(right: 1cm) .. (180: 3.3cm) ;
\draw[thick, blue ,-] (315 : 3.3cm) .. controls +(right: 0cm) and
+(right:0cm) .. (157.5: 3.3cm) ;

\draw[thick, blue ,-] (-135 : 3.3cm) .. controls +(up: 3cm) and
+(right: .1cm)  .. (-157.5: 3.3cm) ;
\draw[thick, blue ,-] (-90 : 3.3cm) .. controls +(up: 3cm) and
+(right: 0cm) .. (-112.5: 3.3cm) ;
\end{scope}

\begin{scope}
\def \n {16}
\def \radius {4cm}
\node[blue] at (0:\radius) {$1$};
\node[blue] at (-360/\n:\radius) {$2$};
\node[blue] at (-360/\n*2:\radius) {$3$};
\node[blue] at (-360/\n*3:\radius) {$4$};
\node[blue] at (-360/\n*4:\radius) {$5$};
\node[blue] at (-360/\n*5:\radius) {$6$};
\node[blue] at (-360/\n*6:\radius) {$7$};
\node[blue] at (-360/\n*7:\radius) {$8$};
\node[blue] at (-360/\n*8:\radius) {$9$};
\node[blue] at (-360/\n*9:\radius) {$10$};
\node[blue] at (-360/\n*10:\radius) {$11$};
\node[blue] at (-360/\n*11:\radius) {$12$};
\node[blue] at (-360/\n*12:\radius) {$13$};
\node[blue] at (-360/\n*13:\radius) {$14$};
\node[blue] at (-360/\n*14:\radius) {$15$};
\node[blue] at (-360/\n*15:\radius) {$16$};
\end{scope}

\begin{scope}
\def \n {16}
\def \radius {2.9cm}
\node[black] at (0:\radius) {$3$};
\node[black] at (-360/\n:\radius) {$4$};
\node[black] at (-360/\n*2:\radius) {$12$};
\node[black] at (-360/\n*3:\radius) {$1$};
\node[black] at (-360/\n*4:\radius) {$7$};
\node[black] at (-360/\n*5:\radius) {$8$};
\node[black] at (-360/\n*6:\radius) {$10$};
\node[black] at (-360/\n*7:\radius) {$9$};
\node[black] at (-360/\n*8:\radius) {$15$};
\node[black] at (-360/\n*9:\radius) {$11$};
\node[black] at (-360/\n*10:\radius) {$16$};
\node[black] at (-360/\n*11:\radius) {$6$};
\node[black] at (-360/\n*12:\radius) {$13$};
\node[black] at (-360/\n*13:\radius) {$14$};
\node[black] at (-360/\n*14:\radius) {$5$};
\node[black] at (-360/\n*15:\radius) {$2$};
\end{scope}

\end{scope}

\begin{scope}[xshift= 10cm]

\begin{scope}
\def \n {16}
\def \radius {3.3cm}
\def \margin {8} 
\foreach \s in {1,...,\n}
{
 \node[blue] at ({360/\n * (\s - 1)}:\radius) {$\bullet$};
 \draw[blue] ({360/\n * (2*\s - 1)}:\radius) arc ({360/\n * (2*\s -
1)}:{360/\n * (2*\s)}:\radius);
}
\end{scope}

\begin{scope}
\draw[thick, blue ,-] (0 : 3.3cm) .. controls +(left: 0cm) and +(down:
0cm) .. (202.5: 3.3cm) ;
\draw[thick, blue ,-] (22.5: 3.3cm) .. controls +(left: 3cm) and
+(down: 0cm) .. (-67.5 : 3.3cm) ;
\draw[thick, blue ,-] (-22.5: 3.3cm) .. controls +(left: 2cm) and
+(down: 0cm) .. (180 : 3.3cm) ;
\draw[thick, blue ,-] (-45 : 3.3cm) .. controls +(left: 0cm) and +(up:
0cm) .. (157.5: 3.3cm) ;
\draw[thick, blue ,-] (45 : 3.3cm) .. controls +(down: 1cm) and +
(down: 3cm)  .. (112.5: 3.3cm) ;
\draw[thick, blue ,-] (135 : 3.3cm) .. controls +(right:1cm) and +
(right: 1cm)  .. (-135: 3.3cm) ;
\draw[thick, blue ,-] (90 : 3.3cm) .. controls +(down:2cm) and +
(down: 0cm)  .. (-112.5: 3.3cm) ;
\draw[thick, blue ,-] (67.5 : 3.3cm) .. controls +(down:2cm) and +
(down: 0cm)  .. (-90: 3.3cm) ;
\end{scope}

\begin{scope}
\def \n {16}
\def \radius {4cm}
\node[blue] at (0:\radius) {$1$};
\node[blue] at (-360/\n:\radius) {$2$};
\node[blue] at (-360/\n*2:\radius) {$3$};
\node[blue] at (-360/\n*3:\radius) {$4$};
\node[blue] at (-360/\n*4:\radius) {$5$};
\node[blue] at (-360/\n*5:\radius) {$6$};
\node[blue] at (-360/\n*6:\radius) {$7$};
\node[blue] at (-360/\n*7:\radius) {$8$};
\node[blue] at (-360/\n*8:\radius) {$9$};
\node[blue] at (-360/\n*9:\radius) {$10$};
\node[blue] at (-360/\n*10:\radius) {$11$};
\node[blue] at (-360/\n*11:\radius) {$12$};
\node[blue] at (-360/\n*12:\radius) {$13$};
\node[blue] at (-360/\n*13:\radius) {$14$};
\node[blue] at (-360/\n*14:\radius) {$15$};
\node[blue] at (-360/\n*15:\radius) {$16$};
\end{scope}

\begin{scope}
\def \n {16}
\def \radius {2.9cm}
\node[black] at (0:\radius) {$2$};
\node[black] at (-360/\n:\radius) {$9$};
\node[black] at (-360/\n*2:\radius) {$16$};
\node[black] at (-360/\n*3:\radius) {$7$};
\node[black] at (-360/\n*4:\radius) {$12$};
\node[black] at (-360/\n*5:\radius) {$13$};
\node[black] at (-360/\n*6:\radius) {$4$};
\node[black] at (-360/\n*7:\radius) {$1$};
\node[black] at (-360/\n*8:\radius) {$10$};
\node[black] at (-360/\n*9:\radius) {$15$};
\node[black] at (-360/\n*10:\radius) {$3$};
\node[black] at (-360/\n*11:\radius) {$5$};
\node[black] at (-360/\n*12:\radius) {$14$};
\node[black] at (-360/\n*13:\radius) {$11$};
\node[black] at (-360/\n*14:\radius) {$6$};
\node[black] at (-360/\n*15:\radius) {$8$};
\end{scope}

\end{scope}

\end{tikzpicture}
\label{G:graphofsigma}
\caption{$\varGamma(x)$ (left) and $\varGamma(y)$ (middle)
and $\varGamma(xy)$ (right).}
\end{center}
\end{figure}
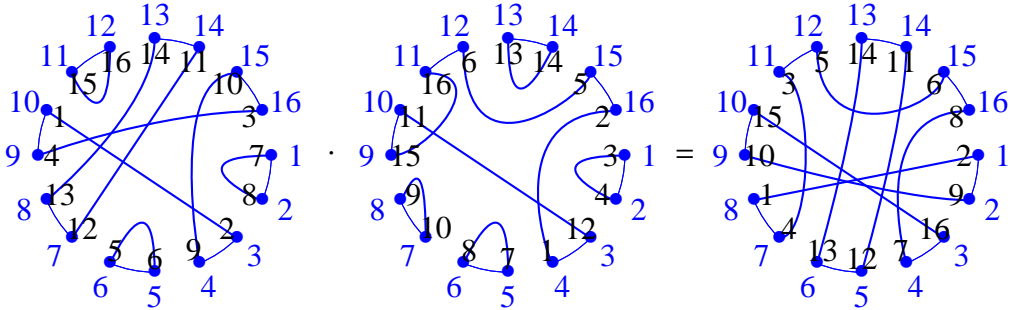

By abuse of terminology, the cycles of the graph $\varGamma(x)$ of an element $x\in S_{2\infty}$ are
called the {\em cycles of $x$}.
For $x$ as in (\ref{A:cycles of x and y}), there are three $1$-cycles, a single $2$-cycle and a single $3$-cycle. 
Note that aside from four $1$-cycles, the element $y$ 
consists of a single $4$-cycle, which we call the \emph{big} cycle of $y$. 
The cycles of the product $xy$ are formed from the the cycles of $x$  by the help of the big cycle of $y$. 
We see that some cycles of $x$ unite together to form a longer cycle of $xy$. 

The $1$-cycles of $y$ do not alter the structure of the cycles of $x$, except moving them around. 
For example, the $1$-cycles $\{5,6,7,8\}$ and $\{ 13,14, 14, 13 \}$ of $y$ transformed the (unique) $2$-cycle of $x$  into the unique $2$-cycle of $xy$.
(See the $2$-cycles in $x$ and $xy$.)
	
Let us understand the effect of the big cycle of $y$ on the cycles of $x$. The interior labels $\{1,2,5,6,16,15,11,12\}$ appear in the big cycle of $y$. 
Images of these numbers under $x$ lie in the same cycle in $xy$. In fact the cycles of $x$ which contains the pairs $\{1,2\}$, 
$\{5,6\}$, $\{11,12\}$, $\{15,16\}$ merged into a single $6$-cycle in $xy$. 
Note that $\{1,2\}$,$\{5,6\}$,  and $\{11,12\}$ belong to different $1$-cycles in $x$, 
and $\{15,16\}$ belongs to a $3$-cycle in $x$. Therefore, the resulting cycle in $xy$ is of size $1+1+1+3 = 6$.

Consequently, four different cycles (three $1$-cycles and a single $3$-cycle) of $x$ merge together and form a single $6$-cycle in $xy$. 
The rest of the cycles of $x$, despite changing their labels, are transformed into cycles of $xy$ preserving their shapes.
These constitute all the cycles of $xy$. 
This is the observation behind the proof of the first part of Theorem \ref{thm:ProductExpnForSingleCycle}.

Let $\tau \in S_{2n}$ be a permutation and let $\alpha$ be a cycle in the graph $\varGamma(\tau)$. 
We introduce a new notation for the cycles of the graph $\varGamma(\tau)$.

Suppose the cycle $\alpha$ is of length $r+1$. Then the cycle $\alpha$ is denoted a sequence of $4(r+1)$ labels, where pairs separated by a 
colon ``:''. 
The sequence starts with an exterior pair, followed by an interior pair, continuing in this fashion alternating between exterior pairs 
and interior pairs until all the vertices of $\alpha$ are visited, traversing the cycle $\alpha$ in clockwise fashion. 
The final pair for $\alpha$ is an always interior pair.
In addition, we call a label $j$ irrelevant to $\alpha$ if $j$ does not belong to any exterior or interior pair of $\alpha$. 
In this new notation, the $3$-cycle $\alpha$ of $y \in K_{(3)}$ of the previous example is written as
$$
\alpha = (3, 4:1, 2:16, 15: 5, 6:12, 11:16, 15: 9, 10:11, 12).
$$
For $\alpha$, the interior labels are $\{ 1, 2, 5, 6, 16, 15, 11, 12 \}$, and $\{ 3,4,16,15,12,11,9,10 \}$ gives the the exterior labels. 
The labels $\{ 7,8,13,14 \}$ are irrelevant to $\alpha$.

\begin{proof} [Proof of the first claim of Theorem \ref{thm:ProductExpnForSingleCycle}.]

Let $n$ be the smallest positive integer such that $y \in K_{(r)} \cap S_{2n}$ and $x \in K_{\lambda} \cap S_{2n}$.
Let $\alpha = (i_1, i_2: i_3, i_4 : \ldots : i_{4r+3}, i_{4r+4} )$ denote the big cycle of $y$. Let $I_\alpha$ be the set of all interior labels of $\alpha$.
Then
\begin{eqnarray*}
xy (i_{4k+2} ) & = & x(i_{4k+3}) \qquad   \text{for}\  k=0,\ldots, r, \label{A:xy.1} \\
xy (i_{4k+1} ) & = & x(i_{4k}) \qquad    \text{for}\   k=1,\ldots, r+1, \label{A:xy.2}\\
xy (i_{1} ) & =  & x(i_{4r+4}). \label{A:xy.3} 
\end{eqnarray*}

Let $\mathcal{A} = \{ A_1,\ldots, A_s \}$ be the set of all cycles of $xy$ determined by the equations (\ref{A:xy.1})--(\ref{A:xy.3}). 
We claim that $\mathcal{A}$ consists of a single cycle, that is $s=1$.

Denote the cycles of $x$ whose set of exterior labels intersects $I_\alpha$ nontrivially by $B_j$ for $j=1, \ldots, t$. 
Clearly, $t \leq r+1$. Denote the lengths of cycles $A_i$ and $B_j$ with $a_i$ and $b_j$ respectively, for $i=1,\ldots, s$ and $j=1,\ldots, t$.

Observe that the union of the sets of interior labels of the cycles $A_i$ for $i=1,\ldots, s$ is precisely the union of the sets of the 
interior labels of $B_j$ for $j=1,\ldots, t$. Therefore,
\begin{align}\label{A:temel1}
\sum^s_{i=1} a_i = \sum^{t}_{j=1}b_j.
\end{align}
Each cycle $C$ of $x$ whose set of exterior labels is disjoint from the set of interior labels of the big cycle $\alpha$ of $y$ 
is transformed into a cycle $C'$ of the same length in $xy$.

In the example above, the cycle $C=(7,8:13,14:13,14:11,12)$ of $x$ is transformed into  the cycle $C' = (5,6:13,14:13,14:11,12)$ of $xy$. 
Enumerate the rest of the cycles of $xy$ and $x$ by $A_{s+i}$ and $B_{t+i}$ respectively, where $B_{t+i}$ of $x$ transforms to 
$A_{s+i}$ of $xy$ for $i=1,\ldots, u$.
Therefore, 
\begin{align*}
\sum^u_{i=1} a_{s+i} = \sum^{u}_{j=1}b_{t+j}.
\end{align*}	

Since $\lgth{ \lambda } + r = \lgth{ \nu }$, we have 
\begin{equation*}
r = \lgth{\nu} - \lgth{\lambda} = \sum^{s+u}_{i=1} (a_i-1) - 
\sum^{t+u}_{j=1}(b_j-1) = \sum^s_{i=1} (a_i-1) - \sum^{t}_{j=1}(b_j-1),
\end{equation*}
which implies $r=t-s$ by Equation \ref{A:temel1}. Since $t \leq r+1$, we conclude that $t = r+1$ and $s =1$. 
This shows that the cycles $B_1, \ldots, B_t$ of $x$ are glued together by the big cycle of $\alpha$ of $y$ and form a single cycle $A_1$ of $xy$.

Since $s=1$, Equation \ref{A:temel1} now reads as $a_1 = \sum_{i =1}^{t} b_i $. Thus 
\begin{itemize}
\item after reordering,  the sequence $(b_i -1)_{i=1}^{t+u}$ is  the partition $\lambda$,
\item after reordering, the sequence
$(a_1 -1 , b_{t+1} -1, \ldots, b_{t+u} -1 )$ is  the partition $\nu$. 
\end{itemize}

Let $\rho \subseteq \lambda$ be the partition corresponding the sequence 
$(b_1-1,\ldots, b_{t}-1)$. Then $ \nu = ( r+\lgth{\rho} )  \cup \lambda - \rho $ as claimed. 
\end{proof}

Inspired by the previous proposition, we introduce the following formalism. Let $\alpha$ be the big cycle of $y \in K_{(r)}$,
and let $B_1,\ldots, B_{t}$ denote the cycles of $x \in K_{\lambda}$ which are blended together to form a single cycle $A_1$ in $K_{\nu}$ by the cycle $\alpha$.
We write 
\begin{align}\label{A:formalism}
B_1 \cdots B_{t} \cdot \alpha = A_1. 
\end{align}

Now fix $z  \in K_{\nu}$. The coefficient of $K_{( r+\lgth{\rho} )  \cup \lambda - \rho} $ in the product expansion of $K_{\lambda} K_{(r)}$ is the number of 
different ways of choosing $x \in K_{\lambda}$ and $y \in K_{(r)}$ so that $z = xy$. 
This coefficient counts the number of ways of choosing the cycle $A_1$ in $z$ and writing it as a product of cycles $B_1, \ldots, B_{t}$ of $x$ 
and the big cycle $\alpha$ of $y$ as it was done in (\ref{A:formalism}).

The number of ways a cycle of length $r+\lgth{ \rho }+1$ can be chosen from $z$ equals $m_{r+\lgth{ \rho }}(\lambda)+1$. 
Notice that after relabeling the entries of the $B_i$'s and $\alpha$ if necessary,  the number of decompositions in (\ref{A:formalism}) 
equals the coefficient of $K_{(r+ \lgth{ \rho })}$ in $K_{\rho} K_{(r)}$. 
By Theorem (\ref{T:sabit}), this coefficient is a constant and is the coefficient of  $K_{(n)}(n+1)$ in $K_{\rho}(n+1) K_{(r)}(n+1)$, where $n= r + \lgth{ \rho }$.

For a partition $\mu$, define $z_\mu := \prod_{i\geq 1} i^{m_i(\mu)} m_i(\mu)!.$
We need to calculate $| K_{\rho}(n) |$ for  a partition $\rho$ such that  $ \mt{wt}(\rho) \leq n$. 
Let $\rho(n)= \rho + (1^{n - \lgth{\rho}})$ be the $n$-completion of $\rho$. 
Then $\rho(n)$ is a partition of $n$ of length $ n - \lgth{\rho}$, for which 
$z_{\rho(n)} =  \prod_{i\geq 0} (i+1)^{m_{i}(\rho)} m_{i} (\rho)!$, where $m_0(\rho) := n - \mt{wt}(\rho)$. Therefore, 
\begin{align*}
| K_{\rho}(n) | &= \frac{| B_n |^2} { 2^{\ell(\rho(n))} z_{\rho(n)} } = \frac{ |B_n|^2} { 2^{ n - \lgth{ \rho }} z_{\rho(n)} }
=  \frac{2^{n+\lgth{\rho}} (n!)^2} { \prod_{i\geq 0} (i+1)^{m_{i}(\rho)} m_{i} (\rho)! }.
\end{align*}
To find an $(r+1)$-cycle $y \in K_{(r)}(n+1)$  such that $w y \in K_{(n)}(n+1)$ we have to choose exterior label pairs from $r+1$ cycles of $w$ 
and then pick one of the $2^{r} r!$ cyclic arrangements of these pairs. 
If $\rho = (1^{m_1(\rho)},2^{m_2(\rho)},\ldots)$, then the number of such $y$ is 
equal to $2^{m_1(\rho)}3^{m_2(\rho)} \cdots 2^{r} r!$.

Finally, the coefficient is the product of $|K_{\rho}(n+1)|$ and $2^{m_1(\rho)}3^{m_2(\rho)} \cdots 2^{r} r!$ divided by
number of elements in $K_{(n)}(n+1)$. Therefore,
\begin{align}\label{A:crucial_step}
\frac{ |K_{\rho}(n+1) | 2^{m_1(\rho)}3^{m_2(\rho)} \cdots 2^{r} r!}{ |K_{(n)}(n+1)| }.
\end{align}

\begin{proof} [Proof of the second claim of Theorem \ref{thm:ProductExpnForSingleCycle}.]
By (\ref{A:crucial_step}), we know that 
\begin{align*}
b_{\lambda \, (r)}^{( r+\lgth{\rho} )  \cup \lambda - \rho} = 
\frac{ (m_{r+\lgth{ \rho }}(\lambda)+1) 2^{m_1(\rho)}3^{m_2(\rho)} \cdots 2^{r} r! |K_{\rho}(n+1) |}{ |K_{(n)}(n+1)| },
\end{align*}
where $n = r+ \lgth{ \rho }$. Now $m_0(\rho) := n+1-\lgth{\rho}-\ell(\rho) = r+1 - \ell(\rho)$. Since		
$| K_{\rho}(n) | = \frac{2^{n+\lgth{\rho}} (n!)^2}{ \prod_{i\geq 0} (i+1)^{m_{i}(\rho)}m_{i} (\rho)! }$, we have

\begin{align*}
b_{\lambda \, (r)}^{( r+\lgth{\rho} )  \cup \lambda - \rho} &=  
\frac{ [  (m_{r+\lgth{ \rho }}(\lambda)+1)\prod_{i\geq 0}  (i+1)^{m_{i}(\rho)}  2^{r} r!   ] 2^{n+1+\lgth{\rho}}}
{  \prod_{i\geq 0}   (i+1)^{m_{i}(\rho)} m_i (\rho)!} \frac{ (n+1) }{ 2^{n+1+n}}\\
&=  \frac{ (m_{r+\lgth{ \rho }}(\lambda)+1)  (n+1)   r!  2^{r+\lgth{\rho}-n}   }{  \prod_{i\geq 0} m_i (\rho)! } \\
&=  \frac{(m_{r+\lgth{ \rho }}(\lambda)+1)(r+\lgth{\rho}+1 ) r!  2^{r+\lgth{\rho}-n}  }{  \prod_{i\geq 0} m_i (\rho)! }\\
&= \frac{(m_{r+\lgth{ \rho }}(\lambda)+1)(r+\lgth{\rho}+1 ) r! }{  \prod_{i\geq 0} m_i (\rho)! }, 
\end{align*}
and the proof is finished. 
\end{proof}

\bibliography{FH.bib}
\bibliographystyle{plain}

\end{document}